\documentclass[a4paper, 10pt]{article}
\usepackage{geometry}
\usepackage[francais, english]{babel}
\usepackage{latexsym,amsfonts, amssymb}
\usepackage{amsthm}
\usepackage[latin1]{inputenc}
\usepackage{amsmath}
\usepackage{amssymb}
\usepackage{amsfonts}
\usepackage{dsfont}
\usepackage{color}
\usepackage{multicol}
\usepackage{hyperref}
\usepackage{amsthm}
\usepackage{latexsym}
\usepackage{fancyhdr}
\usepackage{epsfig}
\usepackage{mathrsfs}
\usepackage{multirow}

\makeatletter
\renewcommand\theequation{\thesection.\arabic{equation}}
\@addtoreset{equation}{section} \makeatother

\newfont{\bbf}{cmbx12 scaled 1435}

\newtheorem{thm}{Theorem}[section]

\newtheorem{lem}{Lemma}[section]

\renewcommand{\theequation}{\thesection.\arabic{equation}}
\newcommand{\eop}{\hspace*{\fill} \ensuremath{\Box}}

\begin{document}

\newcommand{\I}[1]{\mathds{1}_{{#1}}}

\def\Sum{ \displaystyle \sum }
\def\Frac{\displaystyle \frac}

\def\Cit{\mathbb{C}}
\def\esp{\mathbb{E}}
\def\Var{\hbox{\rm Var}}
\def\Cov{\hbox{\rm Cov}}
\def\Supp{\hbox{\rm Supp}}
\def\Card{\hbox{\rm Card}}
\def\Det{\hbox{\rm Det}}
\def\Tr{\hbox{\rm Tr}}
\def\Dim{\hbox{\rm dim}}
\def\Rank{\hbox{\rm dim}}
\def\Id{\hbox{\rm Id}}
\def\Ker{\hbox{\rm Ker}}
\def\ind{\mathbb{I}}
\def\Nit{\mathbb{N}}
\def\Rit{\mathbb{R}}
\def\Zit{\mathbb{Z}}
\def\prob{\mathbb{P}}
\def\where{\rm where}
\def\with{\rm with}
\def\I{\rm I}
\def\J{\rm J}
\def\Ip{\rm I^{\prime}}
\def\Jp{\rm J^{\prime}}
\def\as{\rm a.s}

\setcounter{page}{1} \setlength{\baselineskip}{.32in}

\begin{center}
 {\Large\bf
Nonparametric kernel estimation of the probability density
function of  regression errors using estimated residuals}
\end{center}

%\title{\bf
%Nonparametric kernel estimation of the probability density
%function of  regression errors}

%\begin{center}
%Rawane SAMB
%\end{center}
%\\
%Universit\'e Catholique de Louvain
%\thanks{This article was started
%and finished when I was at Laboratoire de Statistique Théorique et
%Appliquée, Université Pierre et Marie Curie, which support is
%really acknowledged. The author is particularly grateful to
%Emmanuel Guerre and acknowledges him for many helpful comments and
%suggestions. All errors are mine and under my own responsibility.
%}
% \\[.1cm]
% }
%
%\date{\today}

\begin{center}
This version: August 2011
\end{center}

%\maketitle

\begin{abstract}
Consider the nonparametric regression model $Y=m(X)+\varepsilon$,
where the function $m$ is  smooth but unknown, and $\varepsilon$
is independent of $X$. An estimator of the density of the error
term $\varepsilon$ is proposed and its weak consistency is
obtained. The contribution of this paper is twofold. First, we
evaluate the impact of the estimation of the regression function
on the error density estimator. Secondly, the optimal choices of
the first and second step bandwidths used for estimating the
regression function and the error density are proposed. Further,
we investigate  the  asymptotic normality of the error density
estimator and evaluate its performances in simulated examples.

\end{abstract}
\renewcommand{\thefootnote}{\arabic{footnote}}
\setcounter{footnote}{1}
\setlength{\baselineskip}{.26in}
\section{Introduction}

Let $(X_{1},Y_{1}),\ldots,(X_{n},Y_{n})$ be a sample  of
independent replicates of the random vector $(X, Y)$, where $Y$ is
the univariate dependent variable and $X$ is the covariate of
dimension $d$. Let $m(\cdot)$ be the conditional expectation of
$Y$ given $X$ and let
 $\varepsilon$ be
the related regression error term, so that the regression error
model is
\begin{eqnarray}
Y
 =
 m(X)+\varepsilon,
\label{Rm}
\end{eqnarray}
 where $\varepsilon$ is assumed to have mean zero and to be
 statistically independent of $X$, and the function $m(\cdot)$
 is smooth but unknown.
 In this paper, we investigate the problem of nonparametric estimation
of the probability density function (p.d.f) of the error term
$\varepsilon$. The difficulty of  this study is the fact that the
regression error term is not observed  and must  be estimated.
 In such setting, it would be unwise to estimate the error
 density by means of the conditional approach which is based
 on the probability distribution function of  the response variable
 given the covariate. Indeed, this approach is affected by the curse
 of dimensionality, so that the resulting estimator of the residual
  term  would  have considerably a slow rate of convergence if
  the dimension of the explanatory variable is very high.
The strategy used here is based on the estimated residuals, which
are built from the nonparametric estimator of the regression
function $m(\cdot)$. The proposed estimator for the density of
$\varepsilon$ is built by using the estimated residuals as if they
were the true errors, and the weak consistency of this estimator
is obtained. Our results may have many
 possible applications.
 First, the estimator of the density $f(\cdot)$
  of the residual term $\varepsilon$ is an
 important tool for understanding the residuals behavior and
 therefore the fit of the regression model (\ref{Rm}). Indeed,
this estimator  can be used for
 goodness-of-fit tests of a specified error  distribution in a
 parametric or nonparametric regression setting. Some examples can be found in
  Loynes (1980), Akritas and Van  Keilegom (2001),  Cheng and Sun
  (2008). Secondly, the  estimation of $f(\cdot)$ can
   be useful for testing the symmetry of the residuals distribution.
  See Ahmad and Li (1997), Dette, Kusi-Appiah and Neumeyer
 (2002), Neumeyer and Dette (2007) and references therein.
Note also that the estimation of the error density is useful for
forecasting $Y$ by means of a mode approach, since the mode of the
p.d.f of $Y$ given $X$ is $m(x)+\arg\max_{e\in\Rit} f(e)$.
   Another interest in estimating $f(\cdot)$
 is the construction of nonparametric estimators
 of the hazard function of $Y$ given $X$ (see Van Keilegom and Veraverbeke,
 2002), or the estimation of the
density of the response variable $Y$ (see  Escanciano and
 Jacho-Chavez, 2010).

 Many estimators of the p.d.f. of the regression error
$\varepsilon$ can be obtained from estimation of the regression
function and the conditional p.d.f of $Y$ given $X$.  For the
estimation of the latter, see Roussas (1967, 1991) and Youndjé
(1996), among others.  More direct approaches have also been
proposed. Akritas and Van Keilegom (2001)  estimate the cumulative
distribution function of the regression error in heteroscedastic
model with univariate covariates. The estimator they propose is
based on a nonparametric estimation of the residuals. Their
results show the impact of the estimation of the residuals on the
limit distribution of the underlying estimator of the cumulative
distribution function. The results obtained by Akritas and Van
Keilegom (2001) are generalized by Neumeyer and Van Keilegom
(2010) in the case of the same model with multivariate covariates.
M\"{u}ller, Schick and Wefelmeyer (2004) consider the estimation
of moments of the regression error. Quite surprisingly, under
appropriate conditions, the estimator based on the true errors is
less efficient than the  estimator which uses the nonparametric
estimated residuals. The reason is that the latter estimator
better uses the fact that the regression error $\varepsilon$ has
mean zero.  Fu and Yang (2008) study the asymptotic normality of
kernel error density estimators  in parametric nonlinear
autoregressive models. They show that at a fixed point, the
distribution of these error density estimators is normal without
knowing the nonlinear autoregressive function. Wang, Brown, Cai
and Levine (2008) investigate the impact of the estimation of the
regression function on the estimator of the variance function in a
heteroscedastic model. In their study, they show that for a good
estimation of the variance function, it is important to use a very
small bandwidth, and so a weakly biased estimator for the
regression function of their model. Cheng (2005)  establishes the
asymptotic normality of an estimator of $f(\cdot)$ based on the
estimated residuals. This estimator is constructed by splitting
the sample into two parts: the first part is used for the
construction of the estimator of $f(\cdot)$,
 while the second part of the sample is used for the estimation of the
 residuals.  Efromovich (2005) proposes adaptive estimator  of
the error density, based on a density estimator proposed by
Pinsker (1980). Although these authors  used the estimated
residuals for constructing an estimator of the error density, none
of them investigated  the impact of the dimension of the covariate
on the estimation of $f(\cdot)$, nor the influence of the
first-step bandwidth used to estimate $m(\cdot)$, on the estimator
of the error density.

  The contribution of this paper is twofold. First, we
evaluate the impact of the estimation of the regression function
on the error density estimator. Second, the optimal choices of the
first and second step bandwidths used for estimating the
regression function and the residual density respectively, are
proposed. To this end, the difference between the feasible
estimator which uses the estimated residuals, and the unfeasible
one based on the true errors is established. Further,
  we  investigate  the asymptotic normality  of the feasible estimator and
  evaluate its performance through a simulation study.

 The rest of this paper is organized as follows.
Section 2 presents our estimators
 and some notations used in the sequel.  Sections 3 and 4 group our
assumptions and  main results respectively. Section 5 is devoted
to the simulations. Some concluding remarks are given in Section
6, while the proofs of our results are gathered in Section 7 and
in an appendix.

\setcounter{subsection}{0} \setcounter{equation}{0}
\renewcommand{\theequation}{\thesection.\arabic{equation}}

\section{Construction of the estimators and notations}

The approach proposed here for the nonparametric kernel estimation
of $f(e)$ is based on a two-steps procedure, which  builds,
 in a first step, the estimated residuals
\begin{equation}
\widehat{\varepsilon}_i
=
 Y_i-\widehat{m}_{in},
\quad
i=1,\ldots,n,
\label{epschap}
\end{equation}
where $\widehat{m}_{in} =\widehat{m}_{in}(X_i)$ is the leave-one
out
 version of  the Nadaraya-Watson (1964) kernel estimator of
 $m(X_i)$,
\begin{equation}
\widehat{m}_{in} =
 \frac{\sum_{j=1\atop j\neq i}^nY_j
K_0\left(\frac{X_j-X_i}{b_0}\right)}
{\sum_{j=1\atop j\neq i}^n
K_0\left(\frac{X_j-X_i}{b_0}\right)}.
 \label{mchapi}
\end{equation}
Here $K_0(\cdot)$ is a kernel function defined on $\Rit^d$ and
$b_0=b_0(n)$ is a bandwidth sequence. It is tempting to use, in
the second step, the estimated
 $\widehat{\varepsilon}_i$ as if they were the true residuals $\varepsilon_i$.
  This would ignore the fact that the $\widehat{m}_{in}(X_i)$'s
  can result in  severely
 biased estimates of the $m(X_i)$'s for those $X_i$ which are close to the
 boundaries of the support $\mathcal{X}$ of the covariate distribution.
 That is why  our proposed estimator
 trims the observations $X_i$
 outside an inner subset $\mathcal{X}_0$ of
$\mathcal{X}$,
\begin{equation}
\widehat{f}_{n}(e)
 =
\frac{1} {b_1\sum_{i=1}^n \mathds{1}
\left(X_i\in\mathcal{X}_0\right)}
\sum_{i=1}^n \mathds{1}
\left(X_i\in \mathcal{X}_0\right)
K_1\left(
\frac{\widehat{\varepsilon}_i-e}{b_1}
\right),
 \label{fnchap}
\end{equation}
where $K_1(\cdot)$ is a univariate kernel function and
$b_1=b_1(n)$ is a bandwidth sequence.
 This estimator is the so-called two-steps kernel estimator
of $f(e)$. In principle, it would be possible to assume that most
of the $X_i$'s fall in
 $\mathcal{X}_0$ when this set is very close to
$\mathcal{X}$. This would give an estimator close to the more
natural kernel estimator $\sum_{i=1}^n
 K\left((\widehat{\varepsilon}_i-e)/b_1\right)/(nb_1)$.
However, in the rest of the paper, a fixed subset $\mathcal{X}_0$
will be considered for the sake of simplicity.

Observe that the two-steps kernel estimator $\widehat{f}_{n}(e)$
is a feasible estimator in the sense that it does not depend on
any unknown
 quantity, as desirable in practice. This contrasts with the unfeasible
 ideal kernel estimator
\begin{equation}
 \widetilde{f}_{n}(e)
 =
\frac{1} {b_1\sum_{i=1}^n \mathds{1}
\left(X_i\in\mathcal{X}_0\right)}
\sum_{i=1}^n \mathds{1}
\left(X_i \in \mathcal{X}_0\right)
 K_1\left(\frac{\varepsilon_i-e}{b_1}\right),
 \label{fn}
\end{equation}
 which depends in particular on the unknown regression error
terms. It
 is however intuitively clear that a proportion of the estimated residuals
 (those with $X_i$ not close to the boundary of $\mathcal{X}$)
 yield a density estimator rivaling the one based on the
 corresponding proportion of the true errors.

In the sequel we will denote by $\varphi^{(k)}$ the $k^{\rm th}$
derivative of any function  $\varphi$ which is $k$ times
differentiable.

\section{Assumptions}

The assumptions we need for the proofs of the main results are
listed below for convenient reference.

 \vskip 0.3cm \noindent
 {$\bf(A_1)$}
 {\it
 The support $\mathcal{X}$ of $X$  is a subset of
 $\Rit^d$, $\mathcal{X}_0$
 has a nonempty interior and the closure of $\mathcal{X}_0$
 is in the interior of $\mathcal{X}$.
}
\medskip
\\{$\bf(A_2)$}
{\it The p.d.f. $g(\cdot)$ of the i.i.d. covariates $X_i$  is
strictly positive over $\mathcal{X}_0$ and has continuous second
order partial derivatives  over $\mathcal{X}$.
}
\medskip
\\{$\bf(A_3)$}
{\it
 The regression function $m(\cdot)$ has continuous second order partial
derivatives  over  $\mathcal{X}$.
}
\medskip
\\{$\bf(A_4)$}
{\it The i.i.d. centered error regression terms
 $\varepsilon_i$'s have finite 6th moments and are independent of the
 covariates $X_i$'s.
}
\\{$\bf(A_5)$}
{\it
 The  probability  density function $f(\cdot)$ of the $\varepsilon_i$'s
 has  bounded continuous second order
 derivatives over  $\Rit$ and satisfies
 $\sup_{e\in\Rit}|h_p^{(k)}(e)|<\infty$, where  $h_p (e) = e^p f(e)$,
 $p\in[0,2]$ and $k\in\{0,1,2\}$.
  }
\\{$\bf(A_6)$}
{\it The kernel function $K_0$  is  symmetric, continuous over
$\Rit^d$ with support  contained in $[-1/2, 1/2]^d$ and satisfies
$\int\! K_0 (z) dz = 1$.
}
\medskip
\\{$\bf(A_7)$}
{\it
 The kernel function $K_{1}$ is symmetric, has a compact support,  is three times
 continuously differentiable over
 $\Rit$, and satisfies $\int\! K_1 (v) dv = 1$,
 $\int\! K_1^{(\ell)} (v) dv = 0$
 for $\ell=1,2,3$, and
 $\int\! v K_1^{(\ell)}(v) dv=0$ for $\ell=2,3$.
 }
\medskip
\\{$\bf(A_{8})$}
{\it
 The bandwidth $b_0$ decreases to $0$ when $n\rightarrow\infty$ and satisfies,
 for $d^*=\sup\{d+2,2d\}$, $nb_0^{d^*}/\ln n\rightarrow\infty$ and
 $\ln(1/b_0)/\ln(\ln n)\rightarrow\infty$
 when $n\rightarrow\infty$.
}
\\
{$\bf(A_{9})$}
 {\it
The bandwidth $b_1$ decreases to $0$ and satisfies
$n^{(d+8)}b_1^{7(d+4)}\rightarrow\infty$ when
$n\rightarrow\infty$. }

\vskip 0.3cm \noindent
 Assumptions ${\rm (A_2)}$, ${\rm (A_3)}$ and  ${\rm (A_5)}$
  impose that all the
functions
 to be estimated nonparametrically have two bounded derivatives.
 Consequently the conditions  $\int\! z K_0 (z) dz = 0$ and
$\int\!v K_1 (v) dv = 0$, as assumed in ${\rm (A_6)}$ and ${\rm
(A_7)}$, represent
 standard conditions ensuring  that the bias of the resulting
 nonparametric estimators (\ref{mchapi})
and (\ref{fn}) are of order $b_0^2$ and $b_1^2$.
 Assumption ${\rm (A_4)}$ states independence between
the regression error terms
 and the covariates, and the existence of the moments of $\varepsilon$
 up to the sixth order. The interest of this assumption is
 to make easier techniques of proofs for the asymptotic expansion
 of the estimator $\widehat{f}_n(e)$.
The differentiability of $K_1$ imposed in Assumption ${\rm (A_7)}$
is more specific to our
 two-steps estimation method. This assumption is used to expand the
 two-steps kernel estimator
$\widehat{f}_{n}(e)$ in (\ref{fnchap}) around the unfeasible one
$\widetilde{f}_{n}(e)$ from
 (\ref{fn}), using the errors estimation
 $\widehat{\varepsilon}_i - \varepsilon_i$ and the derivatives of $K_1$
 up to the third order.
Assumption ${\rm (A_8)}$ is useful for obtaining the uniform
convergence
 of the Nadaraya-Watson estimator of $m$  (see for
 instance Einmahl and Mason, 2005), and also gives a similar consistency
 result for the leave-one-out estimator $\widehat{m}_{in}$ in
 (\ref{mchapi}). Assumption ${\rm (A_9)}$ is needed in the study of
  the difference between the
 feasible estimator $\widehat{f}_{n}(e)$ and the unfeasible
  estimator $\widetilde{f}_{n}(e)$.

\section{Main results}
This section is devoted to our main results. The first result we
give here concerns the pointwise consistency of the nonparametric
kernel estimator $\widehat{f}_{n}$ of the error  density $f$.
Next, the optimal first-step and second-step bandwidths used to
estimate $f$ are proposed. We will finish this section by
establishing the asymptotic normality of the estimator
$\widehat{f}_{n}$.

\subsection{Pointwise weak consistency}

The following  result gives the order of the difference between
the feasible estimator and the theoretical error density for all
$e\in\Rit$.

\begin{thm}
Under ${\rm (A_1)-(A_9)}$, we have, for all $e\in\Rit$, and
 $b_0$ and $b_1$ going to $0$,
$$
\widehat{f}_{n}(e)-f(e)
 =
 O_{\prob}
\biggl(AMSE(b_1)+R_n(b_0, b_1)\biggr)^{1/2},
$$
 where
$$
AMSE(b_1)
 =
 \esp_n
 \left[
\left( \widetilde{f}_{n}(e) -
 f(e)
 \right)^2
 \right]
 =
 O_{\prob}
\left(b_1^4 + \frac{1}{nb_1}\right),
$$
and
\begin{eqnarray*}
 R_n(b_0, b_1)
 =
 b_0^4
 +
 \left[
 \frac{1}{(nb_1^5)^{1/2}}
 +
 \left(\frac{b_0^d}{b_1^3}\right)^{1/2}
 \right]^2
 \left(
 b_0^4
 +
 \frac{1}{nb_0^d}
 \right)^2
 +
 \left[
 \frac{1}{b_1}
 +
 \left(\frac{b_0^d}{b_1^7}\right)^{1/2}
 \right]^2
 \left(b_0^4+\frac{1}{nb_0^d}\right)^3.
\end{eqnarray*}
 \label{thm1}
\end{thm}

\noindent The result of Theorem \ref{thm1} is based on the
evaluation of the difference between  $\widehat{f}_{n}(e)$ and
$\widetilde{f}_{n}(e)$. This evaluation  gives an indication about
the impact of the estimation of the residuals on the nonparametric
estimation of the regression error density. The remainder term
$R_n(b_0,b_1)$ comes from the replacement of the unknown $m(X_i)$
in $\varepsilon_i$ by the estimate $\widehat{m}_{in}(X_i)$.

\subsection{Optimal first-step and second-step bandwidths
for the pointwise weak consistency}

 As shown in the next result, Theorem \ref{Optimbandw1}
gives some guidelines for the choice of the optimal bandwidth
 $b_0$ used in the nonparametric
estimation of the regression errors. As far as we know, the
optimal choice for $b_0$ has not been investigated before in the
nonparametric literature. In what follows, $a_n \asymp b_n$ means
that $ a_n= O (b_n)$ and $b_n = O(a_n)$, i.e. that there is a
constant $C>0$ such that $|a_n|/C \leq |b_n| \leq C |a_n|$ for $n$
large enough.

\begin{thm}
 Assume
${\rm (A_1)-(A_9)}$ and  define
$$
b_0^* = b_0^*(b_1)
 =
\arg\min_{b_0}
 R_n (b_0, b_1),
$$
where the minimization is performed over bandwidth $b_0$
fulfilling ${\rm (A_8)}$. Then,
$$
b_0^*
\asymp
\max
\left\lbrace
\left(\frac{1}{n^2b_1^3}\right)^{\frac{1}{d+4}},
\left(\frac{1}{n^3b_1^7}\right)^{\frac{1}{2d+4}}
\right\rbrace,
$$
and
$$
R_n(b_0^*, b_1)
\asymp
\max
\left\lbrace
\left(\frac{1}{n^2b_1^3}\right)^{\frac{4}{d+4}} ,
\left(\frac{1}{n^3b_1^7}\right)^{\frac{4}{2d+4}}
 \right\rbrace.
$$
\label{Optimbandw1}
\end{thm}
Our next theorem gives the conditions for which the estimator
$\widehat{f}_{n}(e)$ reaches the optimal rate $n^{-2/5}$ when
$b_0$ takes the value $b_0^*$.
 We prove that for
$d\leq 2$, the bandwidth that minimizes the term
$AMSE(b_1)+R_n(b_0^*, b_1)$ has the same order as $n^{-1/5}$,
yielding the optimal order $n^{-2/5}$ for
$\left(AMSE(b_1)+R_n(b_0^*, b_1)\right)^{1/2}$. Note that the
order $n^{-2/5}$ is the optimal rate achieved  by the optimal
kernel estimator of an univariate
 density. See, for instance,  Bosq and Lecoutre (1987), Scott
(1992) or Wand and Jones (1995).

\begin{thm}
Assume ${\rm (A_1)-(A_9)}$  and let
$$
b_1^*
 =
\arg\min_{b_1}
\biggl(AMSE(b_1)+R_n(b_0^*,b_1)\biggr),
$$
where $b_0^*=b_0^*(b_1)$ is defined as in Theorem
\ref{Optimbandw1}.
 Then,
\begin{enumerate}
\item For $d\leq 2$, we have
$$
b_1^*
\asymp
\left(\frac{1}{n}\right)^{\frac{1}{5}}
$$
and
$$
\biggl(
 AMSE(b_1^*)
+
 R_n(b_0^*, b_1^*)
\biggr)^{\frac{1}{2}}
\asymp
\left(\frac{1}{n}\right)^{\frac{2}{5}}.
$$
\item For $d\geq 3$, we have
$$
b_1^*
\asymp
\left(\frac{1}{n}\right)^{\frac{3}{2d+11}}
$$
and
$$
\biggl(
 AMSE(b_1^*)
 +
 R_n(b_0^*, b_1^*)
\biggr)^{\frac{1}{2}}
\asymp
\left(\frac{1}{n}\right)^{\frac{6}{2d+11}}.
$$
\end{enumerate}
\label{Optimbandw2}
\end{thm}

 The results of Theorem \ref{Optimbandw2} show that the rate
$n^{-2/5}$ is reachable if and only if $d\leq 2$. These results
are derived  from Theorem \ref{Optimbandw1}. This latter implies
that if $b_1\asymp n^{-1/5}$, then $b_0^*$ has the same order as
$$
\max
\left\lbrace
\left(\frac{1}{n}\right)^{\frac{7}{5(d+4)}}
,
\left(\frac{1}{n}\right)^{\frac{8}{5(2d+4)}}
\right\rbrace
=
\left(\frac{1}{n}\right)^{\frac{8}{5(2d+4)}}.
$$
 For $d\leq 2$, this order of $b_0^*$ is smaller
than the one of the optimal bandwidth $\widehat{b}_{0}$ obtained
for the nonparametric kernel estimation of $m(\cdot)$. Indeed, it
has been shown in  Nadaraya (1989, Chapter 4) that the optimal
bandwidth $\widehat{b}_{0}$ needed for the kernel estimation
$m(\cdot)$ satisfies $\widehat{b}_{0}\asymp n^{-1/(d+4)}$.
 For $d=1$, the order of $b_0^*$ is
$n^{-(1/5)\times(4/3)}$ which goes to 0 slightly faster than
$n^{-1/5}$, the optimal order of the bandwidth $\widehat{b}_{0}$.
 For $d=2$, the order of $b_0^*$ is
$n^{-1/5}$. Again this order goes to 0 faster than the order
$n^{-1/6}$ of the optimal bandwidth for the nonparametric kernel
estimation of the regression function with two covariates. This
suggests that for $d=1$ and $d=2$,  the ideal bandwidth $b_0$
needed to estimate the residual terms should be very small.
 Such finding
parallels Wang, Brown, Cai and Levine (2008) who show that a
similar result hold when estimating the conditional variance of a
heteroscedastic regression error term. However Wang et {\it al.}
(2008) do not give the order of the optimal bandwidth to be used
for estimating the regression function in their heteroscedastic
setup.
\\
 For $d\geq 3$, we do not achieve the convergence rate
$n^{-2/5}$ for our proposed estimator $\widehat{f}_n(e)$.
 However, we note that $\widehat{b}_{0}$ goes to $0$ slower than $b_0^*$.
This shows that the convergence rate obtained for
$\widehat{f}_{n}(e)$ is better than the optimal rate achieved in
the case of a classical kernel estimator of a multivariate
density.

\vskip 0.1cm All these results prove that  the best estimator
$\widehat{m}_n$ of $m$ needed for estimating $f$ should use a very
small bandwidth $b_0$. This suggests that $\widehat{m}_n$ should
be less biased and should have a higher variance than the optimal
nonparametric kernel estimation of $m$. Consequently the
estimators of $m$ with smaller bias should be preferred in our
framework, compared to the case where the regression function $m$
is the parameter of interest. Indeed, in our case, as in Wang et
{\it al.} (2008), the square of the bias is of more important than
the variance.

\subsection{Asymptotic normality}
Our last result concerns the  asymptotic normality of the
estimator $\widehat{f}_{n}(e)$.

\begin{thm}
Assume ${\rm (A_1)-(A_{9})}$ and
$$
{(\rm\bf{A}_{10}):}
\quad nb_0^{d+4}=O(1),
\quad nb_0^4b_1=o(1),
\quad
 nb_0^{d}b_1^3\rightarrow\infty,
$$
when $n$ goes to $\infty$. Then,
$$
\sqrt{nb_1}
\left(
 \widehat{f}_{n}(e)
 -
\overline{f}_{n}(e) \right) \stackrel{d}{\rightarrow} N\left(
 0,
\frac{f(e)} {\prob\left(X\in\mathcal{X}_0\right)} \int K_1^2(v)
 dv
 \right),
$$
where
$$
\overline{f}_{n}(e)
 =
f(e) + \frac{b_1^2}{2}f^{(2)}(e) \int
 v^2 K_1(v) dv
+
o\left(b_1^2\right).
$$
\label{normalite}
\end{thm}

\noindent Note that for $d\leq 2$, $b_1=b_1^*$ and $b_0=b_0^*$,
Theorems \ref{Optimbandw1} and \ref{Optimbandw2} imply that
$$
b_1
\asymp
\left(\frac{1}{n}\right)^{\frac{1}{5}},
\quad
b_0
\asymp
\left(\frac{1}{n}\right)^{\frac{8}{5(2d+4)}},
$$
which yields
$$
nb_0^{d+4}
 \asymp
\left(\frac{1}{n}\right)^{\frac{12-2d}{5(2d+4)}},
\quad
nb_0^4b_1
\asymp
\left(\frac{1}{n}\right)^{\frac{16-8d}{5(2d+4)}},
 \quad
nb_0^db_1^3
\asymp
\left(\frac{1}{n}\right)^{\frac{4d-8}{5(2d+4)}}.
$$
This shows that for $d=1$,  Assumption $(\rm A_{10})$ is
realizable with the  bandwidths $b_0^*$ and $b_1^*$. But with
these bandwidths, the last constraint of $(\rm A_{10})$ is not
satisfied for $d=2$, since $nb_0^db_1^3$ is bounded as
$n\rightarrow\infty$.

\section{Simulations}
In this section we report simulation results evaluating the finite
sample behavior of the estimators $\widetilde{f}_n$ and
$\widehat{f}_n$. In two examples, we evaluate the performance of
these estimators in terms of asymptotic biases, variances and mean
square errors. The first example concerns a one-dimensional
regression model (univariate covariate), while the second example
is devoted to a regression model with a three-dimensional
covariate.

\subsection{Univariate case}
We work with the following data generating model
\begin{equation}
Y=1+\sin(\pi X)+\varepsilon,
 \label{model1}
\end{equation}
where $\varepsilon\sim N(0,1)$ and $X\sim U[0,1]$. We use the
kernel $K=K_j(x)=(15/16)(1-x^2)^2
 \mathds{1}(|x|\leq 1)$ $(j=0,1)$. Our results are based on  $300$
 simulation runs. For the bandwidth choice, we consider the results
 of Theorems \ref{Optimbandw1} and \ref{Optimbandw2} and take
$$
b_1
=
\widetilde{b}_1,
\quad
b_0
=
c_0
\times
\max
\left\lbrace
\left(\frac{1}{n^2b_1^3}\right)^{\frac{1}{d+4}},
\left(\frac{1}{n^3b_1^7}\right)^{\frac{1}{2d+4}}
\right\rbrace
=
c_0\left(\frac{1}{n^3b_1^7}\right)^{\frac{1}{6}},
$$
where $d=1$, $c_0$ is a given constant in $[0,1]$ and
$\widetilde{b}_1=1.06\times\widetilde{\sigma}_{\varepsilon}\times
n^{-1/5}$ is the Silverman's (1986) rule of thumb bandwidth for
the estimator $\widetilde{f}_n$. Here
$\widetilde{\sigma}_{\varepsilon}$ is the average standard
deviation of the generated errors. For the estimators
$\widetilde{f}_n$ and $\widehat{f}_n$, we  consider
$\mathcal{X}_0=[\delta,1-\delta]$, $\delta=0.001$.

\vskip 0.3cm\noindent {\large{TABLE 1 HERE}}

\vskip 0.3cm  In Table \ref{Table1} we give some values of the
bias, variance and mean square error of $\widehat{f}_n(e)$ at the
points $e=-1$, $0$ and $1$ for different sample sizes. For each
sample, the values are calculated for $c_0=0.25, 0.5$ and $1$.
From Table \ref{Table1} we see that
 our method seems to work
well, since the variance and mean square error of
$\widehat{f}_n(e)$ are very close to $0$. We also observe that the
performance of $\widehat{f}_n(e)$ should not be very sensitive to
the choice of the constant $c_0$, since the variations of the
variance and the mean square error are practically negligible.
Further, we note that for $e=-1, 1$ and $n=100$ the variance and
the mean square error of $\widehat{f}_n(e)$ are smaller than the
ones of $\widetilde{f}_n(e)$. This fact parallels the surprising
situation noticed in Müller, Schick and Wefelmeyer (2004) for the
nonparametric kernel estimation of moments of the regression
error.

\vskip 0.3cm\noindent {\large{FIGURE 1 HERE}}

\begin{figure}[H]
\begin{center}
\includegraphics[width=11cm, height=7.5cm]{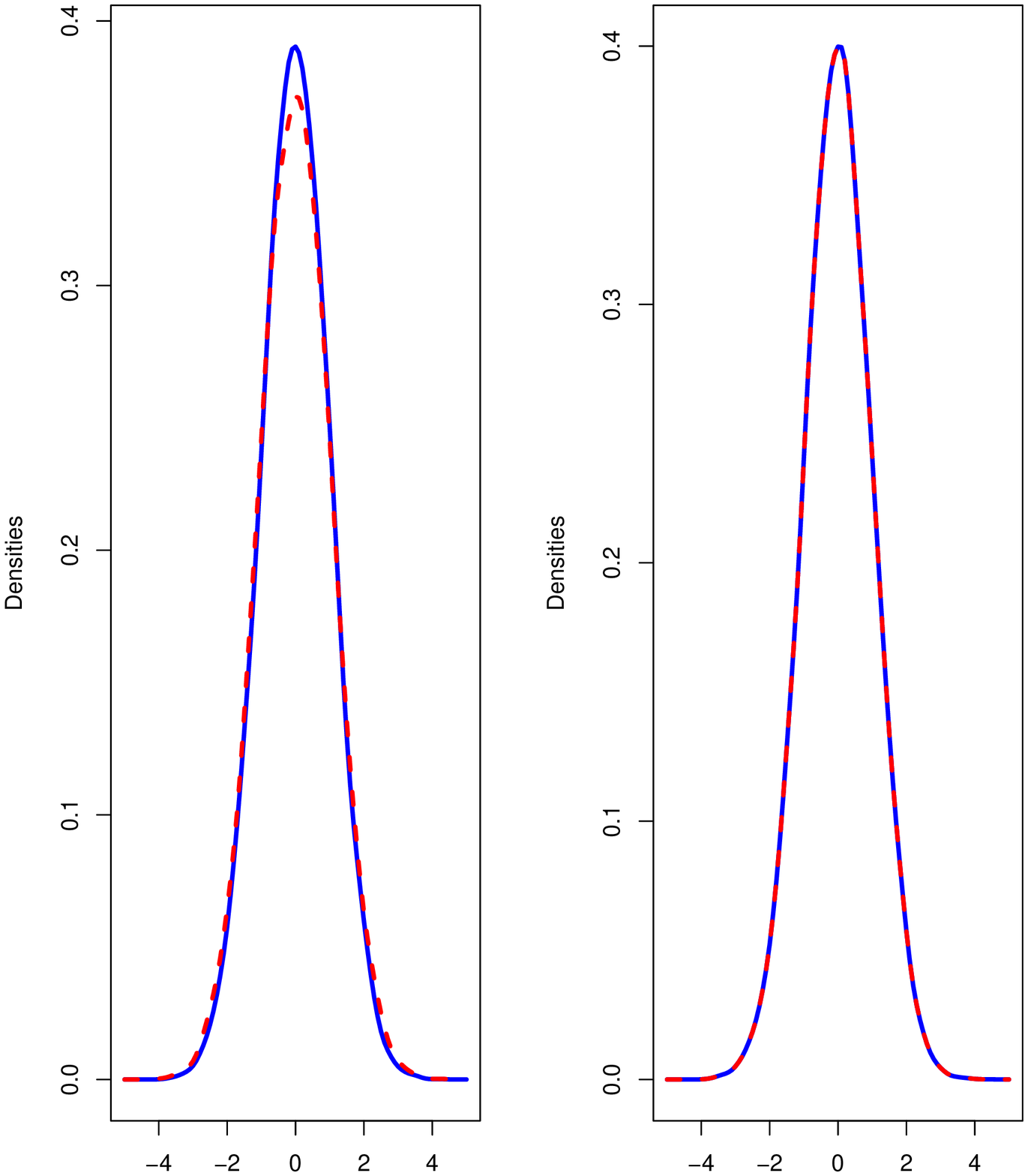}
\end{center}
 \caption{Curves of the densities
$\widetilde{f}_n$ (dashed line) and $\widehat{f}_n$ (solid line)
in univariate case for $c_0=1$ and for sample sizes $n=50$ (left
side) and $n=100$ (right side). All the values of
$\widetilde{f}_n$ and $\widehat{f}_n$ are calculated from $300$
replicates of generated data.}
 \label{Univdensity}
\end{figure}

\vskip 0.3cm Figure \ref{Univdensity} compares the curves of
$\widetilde{f}_n$ and $\widehat{f}_n$ for $c_0=1$ and for samples
size $n=50$ and $n=100$. We observe almost no difference between
the performances of these two estimators. This should suggest that
the estimators $\widetilde{f}_n$ and $\widehat{f}_n$ are
asymptotically equivalent when $n\rightarrow\infty$.

%\begin{figure}[htbp]

%\begin{table}[htbp]
\begin{table}[H]
\begin{center}
\begin{tabular}[b]{|c|c|c|c|c|c|c|c|c|}
\hline $e$ & $n$& $c_0$& \multicolumn{3}{|c|}{ Estimator
$\widetilde{f}_n$}& \multicolumn{3}{|c|}{ Estimator
$\widehat{f}_n$}
\\
\cline{4-9}  &  &  & Bias & Variance & MSE & Bias & Variance & MSE
\\\hline
  $-1$
&\multirow{3}{1cm}{50}
 &$0.25$& 0.2380 & 0.0080 & 0.0647 & 0.1592& 0.0062 & 0.0316
 \\
 \cline{3-9}
 &  & 0.5 & 0.2380  & 0.0080 & 0.0647 &0.2185 &0.0069 & 0.0547
 \\
 \cline{3-9}
 & & 1  &0.2380  & 0.0080 & 0.0647 & 0.2357& 0.0071& 0.0627
   \\
\cline{2-9} &\multirow{3}{1cm}{100}
 &$0.25$&-0.0019  & 0.0034 & 0.0034 &-0.0038 & 0.0027 & 0.0027
 \\
 \cline{3-9}
 &  & 0.5 & -0.0019 & 0.0034 & 0.0034 & -0.0026&0.0034 & 0.0034
 \\
 \cline{3-9}
 & & 1  & -0.0019 & 0.0034 & 0.0034 & 0.0022 & 0.0030 & 0.0030
 \\\hline
$0$
&\multirow{3}{1cm}{50}
 &$0.25$ &0.3843  & 0.0106 & 0.1583 & 0.1291&0.0111 & 0.0278
 \\
 \cline{3-9}
 &  & 0.5 & 0.3843 & 0.0106 & 0.1583 & 0.2391& 0.0079& 0.0646
 \\
 \cline{3-9}
 & & 1  &0.3843 & 0.0106 & 0.1583 & 0.2886& 0.0104& 0.0937
   \\
\cline{2-9} &\multirow{3}{1cm}{100}
 &$0.25$ &0.0008 &0.0054  &  0.0054 &-0.0440 & 0.0044 & 0.0063
 \\
 \cline{3-9}
 &  & 0.5 & 0.0008 & 0.0054 & 0.0054 &-0.0242 & 0.0053 & 0.0059
 \\
 \cline{3-9}
 & & 1  & 0.0008 &0.0054  & 0.0054 &-0.0137 & 0.0050 & 0.0062
 \\\hline
$1$
&\multirow{3}{1cm}{50}
 &$0.25$ &0.2391  &0.0079  & 0.0651 &0.1557 &0.0579 & 0.0300
 \\
 \cline{3-9}
 &  & 0.5 & 0.2391 & 0.0079 & 0.0651 & 0.2122&0.0069 & 0.0520
 \\
 \cline{3-9}
 & & 1  &0.2391  & 0.0079 & 0.0651 & 0.2275& 0.0071& 0.0589
   \\
\cline{2-9} &\multirow{3}{1cm}{100}
 &$0.25$ &-0.0007  & 0.0038 & 0.0038 & -0.0042& 0.0033& 0.0033
 \\
 \cline{3-9}
 &  & 0.5 & -0.0007 & 0.0038 & 0.0038 &-0.0058 & 0.0034& 0.0035
 \\
 \cline{3-9}
 & & 1  & -0.0007 & 0.0038 & 0.0038 & -0.0063& 0.0033 & 0.0034
  \\\hline
\end{tabular}
 \end{center}
 \caption{The table compares some values of the bias, variance and
 mean square error of the estimators $\widetilde{f}_n$ and
 $\widehat{f}_n$ when the data are generated from
 Model \ref{model1}. All these
 values are based on $300$ simulations runs.}
\label{Table1}
\end{table}

\subsection{Trivariate case}
We consider the model
\begin{equation}
Y=1+X_1+X_2^2+\sin(\pi X_3)+\varepsilon,
\label{model2}
\end{equation}
where $\varepsilon\sim N(0,1)$ and $X_1,X_2,X_3\sim U[0,1]$. As in
the univariate case, our study is based  on $300$ simulation runs.
We use the kernels $K_1(x)=(15/16)(1-x^2)^2
 \mathds{1}(|x|\leq 1)$, $K_0(x_1,x_2,x_3)=\prod_{j=1}^3K_1(x_j)$
 and consider $\mathcal{X}_0=[\delta,1-\delta]^3$, $\delta=0.001$.
 We use the bandwidths
$$
b_1
=
\widetilde{b}_1,
\quad
b_0
=
c_0
\times
\max
\left\lbrace
\left(\frac{1}{n^2b_1^3}\right)^{\frac{1}{d+4}},
\left(\frac{1}{n^3b_1^7}\right)^{\frac{1}{2d+4}}
\right\rbrace
=
c_0\left(\frac{1}{n^2b_1^3}\right)^{\frac{1}{7}},
$$
where $d=3$, $c_0\in[0,1]$ and $\widetilde{b}_1$ is the average
standard deviation on the generated errors.

\vskip 0.3cm\noindent {\large{FIGURE 2 HERE}}

\begin{figure}[htbp]
%\begin{figure}[H]
\begin{center}
\includegraphics[width=11cm, height=7.6cm]{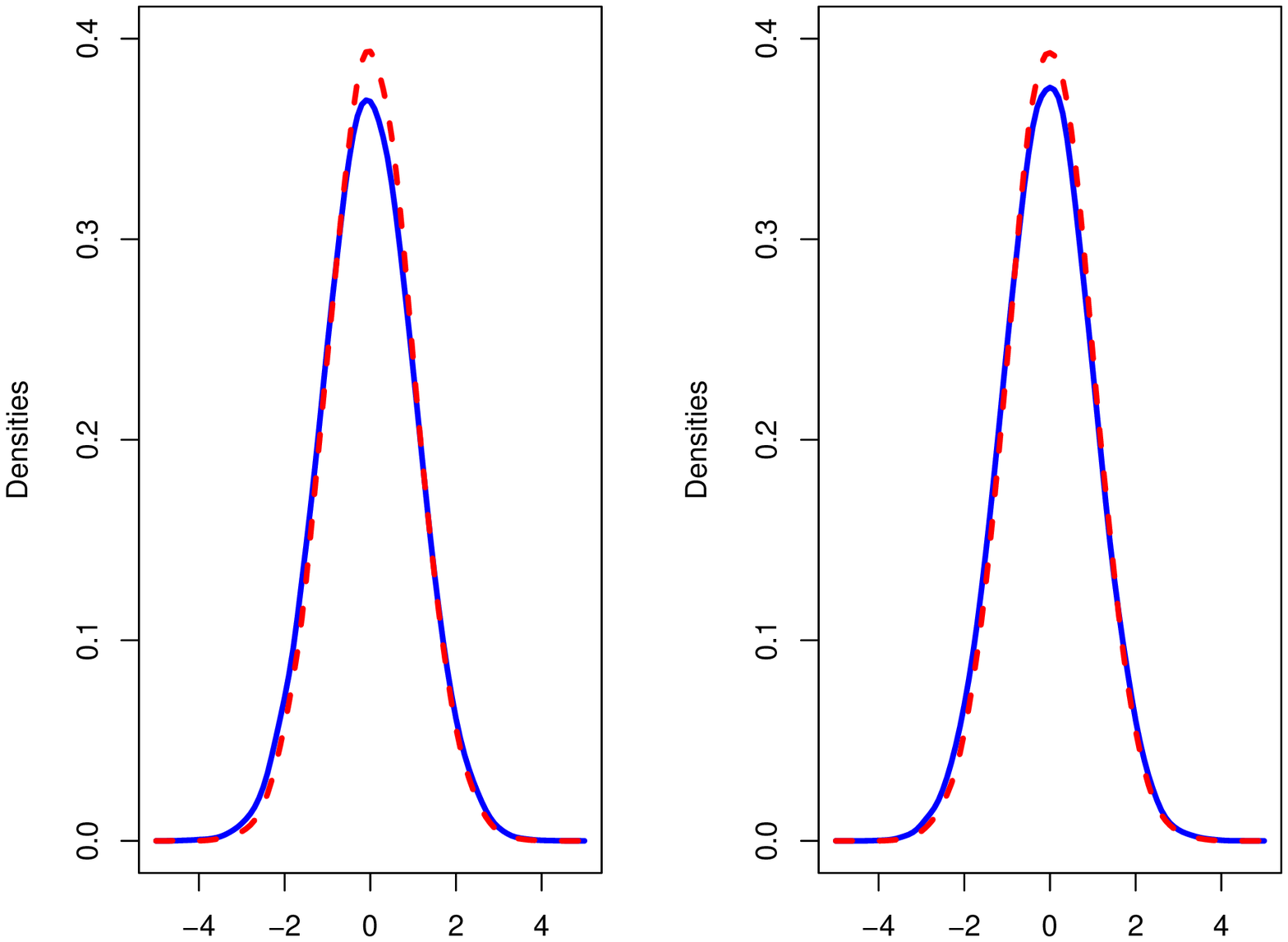}
\end{center}
 \label{Bivdensity}
 \caption{Curves of the densities
$\widetilde{f}_n$ (dashed line) and $\widehat{f}_n$ (solid line)
in trivariate case for $c_0=1$ and for sample sizes $n=100$ (left
side) and $n=200$ (right side). The values are computed from $300$
simulation runs.}
\end{figure}

\vskip 0.3cm Figure \ref{Bivdensity} compares the curves of
$\widetilde{f}_n(e)$ and $\widehat{f}_n(e)$ for $c_0=1$ and sample
sizes $n=100$ and $n=200$. We note a difference between the curves
at the neighborhood of the inflexion point $e=0$. But this
difference is less important for $n=200$. This augurs that for $e$
very close to $0$, the difference between $\widehat{f}_n(e)$ and
$\widetilde{f}_n(e)$ should be negligible only when the size of
the samples is large enough.

\vskip 0.3cm\noindent {\large{TABLE 2 HERE}}

\vskip 0.3cm In Table \ref{Table2} we give some values of the
bias, variance and mean square error of $\widetilde{f}_n(e)$ and
$\widehat{f}_n(e)$ for $c_0=0.25$, $0.5$ and $1$. We see
 that the mean square error of $\widehat{f}_n(e)$ is
greater than the one of $\widetilde{f}_n(e)$.
 Further, we observe that the performance of
$\widehat{f}_n(e)$ should be sensitive to the choice of the
constant $c_0$. For example, for $e=0$, $c_0=0.5$ and $c_0=1$, the
mean square error of $\widehat{f}_n(e)$ is very high compared to
the sum of the variance and the square of the bias.

%\begin{table}[htbp]
\begin{table}[H]
\begin{center}
\begin{tabular}[b]{|c|c|c|c|c|c|c|c|c|}
\hline $e$ & $n$& $c_0$& \multicolumn{3}{|c|}{ Estimator
$\widetilde{f}_n$}& \multicolumn{3}{|c|}{ Estimator
$\widehat{f}_n$}
\\
\cline{4-9}  &  &  & Bias & Variance & MSE & Bias & Variance & MSE
\\\hline
  $-1$
&\multirow{3}{1cm}{100}
 &$0.25$& -0.0013 & 0.0035 &0.0036 & -0.1250 & 0.0015 & 0.0539
 \\
 \cline{3-9}
 &  & 0.5 & -0.0013 & 0.0035 & 0.0036 & -0.0180 & 0.0027& 0.0416
 \\
 \cline{3-9}
 & & 1  & -0.0013 & 0.0035 & 0.0036 & 0.0064 & 0.0027 & 0.0483
   \\
\cline{2-9} &\multirow{3}{1cm}{200}
 &$0.25$& 0.0020 & 0.0019 & 0.0019 &-0.1078 & 0.0011 & 0.0531
 \\
 \cline{3-9}
 &  & 0.5 & 0.0020 & 0.0019 & 0.0019 &-0.0114 & 0.0014& 0.0431
 \\
 \cline{3-9}
 & & 1  & 0.0020 & 0.0019 & 0.0019 & 0.0017 & 0.0015 & 0.0525
 \\\hline
$0$ &\multirow{3}{1cm}{100}
 &$0.25$ & -0.0049 & 0.0047 & 0.0047 & -0.1858 & 0.0042 & 0.0329
 \\
 \cline{3-9}
 &  & 0.5 & -0.0049 & 0.0047 & 0.0047 & -0.0930 & 0.0045 & 0.1451
 \\
 \cline{3-9}
 & & 1  & -0.0049 & 0.0047 & 0.0047 & -0.0377 & 0.0044 & 0.3318
   \\
\cline{2-9} &\multirow{3}{1cm}{200}
 &$0.25$ & -0.0024 & 0.0030 & 0.0030 & -0.1713 & 0.0028 & 0.0378
 \\
 \cline{3-9}
 &  & 0.5 & -0.0024 & 0.0030 & 0.0030 & -0.0764 & 0.0030& 0.1817
 \\
 \cline{3-9}
 & & 1  & -0.0024  & 0.0030 & 0.0030 & -0.0297 & 0.0026 & 0.3591
 \\\hline
$1$ &\multirow{3}{1cm}{100}
 &$0.25$ & -0.0020 & 0.0031 & 0.0031 & 0.0341 &0.0031 & 0.0419
 \\
 \cline{3-9}
 &  & 0.5 & -0.0020 & 0.0031 & 0.0031 &-0.0131 & 0.0025 & 0.0033
 \\
 \cline{3-9}
 & & 1  & -0.0020 & 0.0031 & 0.0031 & -0.0010 & 0.0028 & 0.0416
   \\
\cline{2-9} &\multirow{3}{1cm}{200}
 &$0.25$ & -0.0064 & 0.0019 & 0.0019 &0.0239 & 0.0020 & 0.0325
 \\
 \cline{3-9}
 & & 0.5 &-0.0006 & 0.0019 & 0.0019 & -0.0101 & 0.0016 & 0.0062
 \\
 \cline{3-9}
 & & 1  &-0.0006 & 0.0019 & 0.0019 & -0.0126 & 0.0016 & 0.0477
  \\\hline
\end{tabular}
\end{center}
 \caption{The table gives some values of the bias, variance and
 mean square error of  $\widetilde{f}_n$ and
 $\widehat{f}_n$  when data are generated from
 Model \ref{model2}.
All values are based on $300$ replications of simulated data.}
\label{Table2}
\end{table}

\section{ Conclusion}
 The aim of this paper was to investigate
  the nonparametric kernel estimation of the
probability density function of the regression error using the
estimated residuals. First, we evaluated the impact of the
estimation of the regression function on the error density
estimator. To this aim, the difference between the feasible
density estimator based on the estimated residuals and the
unfeasible one using the true errors was investigated.
 Second, the optimal choices of the first and second
step bandwidths used for estimating the regression function and
the error density were proposed. Further, we establish the
asymptotic normality of the feasible estimator. The strategy used
here to estimate the error density is based on a two-steps
procedure which, in a first step, replaces the unobserved
residuals terms by some nonparametric estimators
$\widehat{\varepsilon}_i=Y_i-\widehat{m}_{n}(X_i)$, where
$\widehat{m}_{n}(X_i)$ is a nonparametric estimator of $m(X_i)$.
 In a second step, the estimated residuals $\widehat{\varepsilon}_i$ are
used to estimate the error density $f(\cdot)$, as if they were the
true $\varepsilon_i$'s. Though  proceeding may remedy the curse of
dimensionality for large sample sizes, a challenging issue was to
evaluate the impact of the estimated residuals on the estimation
of $f(\cdot)$, and to find the order of the optimal first-step
bandwidth $b_0$ used for estimating the error terms. For the
choice of $b_0$, our results show that the ideal bandwidth for
$b_0$ should be smaller than the optimal bandwidth for the
nonparametric kernel estimation of $m(\cdot)$. This means that the
best estimator of $m(\cdot)$ needed for estimating $f(\cdot)$
should have a lower bias and a higher variance than the
 classical kernel regression estimator. With
this ideal choice of $b_0$, we establish that for $d\leq 2$, the
 estimator $\widehat{f}_{n}(e)$ of
$f(e)$ can attain the convergence rate $n^{-2/5}$, which
corresponds to the optimal consistency rate achieved by  the
univariate kernel density estimator. For $d\geq 3$, the rate
$n^{-2/5}$ is not reachable by our estimator $\widehat{f}_n(e)$.
 However, the rate we obtain for
$\widehat{f}_{n}(e)$ is better than the optimal one achieved in
the case of the kernel estimation of a multivariate density.

\setcounter{subsection}{0} \setcounter{equation}{0}
\renewcommand{\theequation}{\thesection.\arabic{equation}}

\section{ Proofs section}

\subsection*{Intermediate Lemmas for Theorem \ref{thm1}}

\begin{lem}
Define, for $x\in\mathcal{X}_0$,
$$
\widehat{g}_n(x)
 =
 \frac{1}{nb_0^d}
 \sum_{i=1}^n
 K_0
 \left(
 \frac{X_i-x}{b_0}
 \right),
 \quad
 \overline{g}_n(x)
 =
 \esp\left[\widehat{g}_n(x)\right].
$$
Then under ${\rm (A_1)-(A_2)}$, ${\rm (A_6)}$ and ${\rm (A_8)}$,
 we have, when $b_0$ goes to $0$,
$$
\sup_{x\in\mathcal{X}_0}
\left|\overline{g}_{n}(x)-g(x)\right|
 =
O\left(b_0^{2}\right),
\quad
\sup_{x\in\mathcal{X}_0}
 \left|
\widehat{g}_{n}(x)
 -
 \overline{g}_n(x)
\right|
 = O_{\prob}
 \left(
 b_0^4 +
 \frac{ \ln n}{nb_0^d}
 \right)^{1/2},
$$
 and
$$
\sup_{x\in\mathcal{X}_0}
 \left|
 \frac{1}{\widehat{g}_{n}(x)}
   -
\frac{1}{g(x)}
\right|
 =
 O_{\prob}
 \left(
 b_0^4
 +
 \frac{\ln n}{nb_0^d}
\right)^{1/2}.
$$
\label{Estig}
\end{lem}

\begin{lem}
Set
$$
f_{in}(e) = \frac{\mathds{1}\left(X_i \in\mathcal{X}_0\right)}
{b_1\prob\left(X\in\mathcal{X}_0\right)}
K_1\left(\frac{\varepsilon_i-e}{b_1}\right).
$$
Then under  ${\rm (A_4)}$, ${\rm (A_5)}$ and ${\rm (A_7)}$, we
have, for $b_1$ going to $0$, and for some constant $C>0$,
\begin{eqnarray*}
\esp
 f_{in}(e)
&=& f(e) + \frac{b_1^2}{2} f^{(2)}(e) \int v^2 K_1(v) dv +
o\left(b_1^2\right),
\\
\Var \left(f_{in}(e)\right) &=& \frac {f(e)}
{b_1\prob\left(X\in\mathcal{X}_0\right)} \int K_1^2(v) dv +
o\left(\frac{1}{b_1}\right),
\\
\esp \left| f_{in}(e) - \esp
 f_{in}(e)
\right|^3
&\leq&
 \frac{
C f(e)}
 {b_1^2\prob^2\left(X\in\mathcal{X}_0\right)}
 \int
 \left|
 K_1(v)
\right|^3
dv
+
 o\left(\frac{1}{b_1^2}\right).
\end{eqnarray*}
\label{Momfin1}
\end{lem}

\begin{lem}
Define
\begin{eqnarray*}
S_n
 &=&
 \sum_{i=1}^n
 \mathds{1}
 \left(X_{i} \in \mathcal{X}_0\right)
 \left(\widehat{m}_{in}-m(X_{i})\right)
 K_1^{(1)}
 \left(\frac{\varepsilon_{i}-e}{b_1}\right),
 \\
 T_n
 &=&
 \sum_{i=1}^n
 \mathds{1}
 \left(X_{i} \in \mathcal{X}_0\right)
 \left(\widehat{m}_{in}-m(X_{i})\right)^2
 K_1^{(2)}
 \left(\frac{\varepsilon_{i}-e}{b_1}\right),
 \\
 R_n
 &=&
 \sum_{i=1}^n
 \mathds{1}
 \left( X_{i}
 \in\mathcal{X}_0\right)
 \left(\widehat{m}_{in}-m(X_{i})\right)^3
 \int_{0}^{1}
 (1-t)^2
 K_1^{(3)}
 \left(
 \frac{
 \varepsilon_i-t(\widehat{m}_{in} -m(X_i))-e}{b_1}
 \right)
 dt.
\end{eqnarray*}
Then under ${\rm (A_1)-(A_9)}$, we have, for $b_0$ and $b_1$ small
enough,
\begin{eqnarray*}
S_n
 &=&
 O_{\prob}
\left[
b_0^2\left(nb_1^2+(nb_1)^{1/2}\right)
 +
\left(nb_1^4+\frac{b_1}{b_0^d}\right)^{1/2}
 \right],
 \\
 T_{n}
 &=&
 O_{\prob}
\left[
\left(
nb_1^3+\left(nb_1\right)^{1/2}
+
\left(n^2b_0^db_1^{3}\right)^{1/2}
\right)
\left(b_0^4+\frac{1}{nb_0^d}\right)
\right],
\\
R_n
&=&
O_{\prob}
\left[
\left(
nb_1^3
+
\left(n^2b_0^db_1\right)^{1/2}
 \right)
\left(b_0^4+\frac{1}{nb_0^d}\right)^{3/2}
\right].
\end{eqnarray*}
\label{STR}
\end{lem}

\begin{lem}
Under  ${\rm (A_5)}$ and ${\rm (A_7)}$ we have, for some constant
$C>0$, and for any $e$ in $\Rit$ and $p\in[0,2]$,
\begin{eqnarray}
\left|
 \int
 K_1^{(1)}
 \left( \frac{e-e}{b_1} \right)^2
 e^pf(e) de
 \right|
 \leq C b_1,
 &&
 \left|
  \int
  K_1^{(1)}
  \left(
 \frac{e-e}{b_1}
 \right)
 e^pf(e)de
  \right|
 \leq C b_1^2,
\label{MomderK1}
\\
\left| \int
 K_1^{(2)}
 \left( \frac{e-e}{b_1} \right)^2
  e^pf(e)
  de
\right|
 \leq
  C b_1,
  &&
 \left|
\int K_1^{(2)}
 \left(
 \frac{e-e}{b_1}
 \right)
 e^pf(e)de
\right|
 \leq C b_1^3,
 \label{MomderK2}
\\
\left| \int
 K_1^{(3)}
 \left( \frac{e-e}{b_1} \right)^2
  e^pf(e)
  de
\right|
 \leq
  C b_1,
  &&
 \left|
\int K_1^{(3)}
 \left(
 \frac{e-e}{b_1}
 \right)
 e^pf(e)de
\right|
 \leq
 C b_1^3.
 \label{MomderK3}
\end{eqnarray}
\label{MomderK}
\end{lem}

\begin{lem}
Let
$$
\beta_{in}
 =
\frac{\mathds{1}\left(X_i\in\mathcal{X}_0\right)}
{nb_0^d\widehat{g}_{in}}
\sum_{j=1, j\neq i}^n
\left(m(X_j)-m(X_i)\right)
K_0\left(\frac{X_j-X_i}{b_0}\right),
$$
where
$$
\widehat{g}_{in}
=
\frac{1}{nb_0^d}
\sum_{j=1\atop j\neq i}^n
K_0\left(\frac{X_j-X_i}{b_0}\right).
$$
Then, under ${\rm (A_1)-(A_9)}$, we have, when $b_0$ and $b_1$ go
to $0$,
\begin{eqnarray*}
\sum_{i=1}^n \beta_{in}
 K_1^{(1)}
 \left(
 \frac{\varepsilon_i-e}{b_1}
 \right)
 =
 O_{\prob}
 \left(b_0^2\right)
 \left(nb_1^2+(nb_1)^{1/2}\right).
\end{eqnarray*}
\label{Betasum}
\end{lem}

\begin{lem}
Let
$$
\Sigma_{in}
 =
\frac{\mathds{1}\left(X_i\in\mathcal{X}_0\right)}
{nb_0^d\widehat{g}_{in}} \sum_{j=1, j\neq i}^n
 \varepsilon_j
K_0\left(\frac{X_j-X_i}{b_0}\right).
$$
Then, under ${\rm (A_1)-(A_9)}$, we have
\begin{eqnarray*}
\sum_{i=1}^n \Sigma_{in}
 K_1^{(1)}
 \left(
 \frac{\varepsilon_i-e}{b_1}
 \right)
 =
 O_{\prob}
 \left(
 nb_1^4
 +
 \frac{b_1}{b_0^d}
\right)^{1/2}.
\end{eqnarray*}
\label{Sigsum}
\end{lem}

\begin{lem}
Let $\esp_n[\cdot]$ be the conditional mean given
$X_1,\ldots,X_n$. Then under ${\rm (A_1)-(A_{9})}$,
 we have
\begin{eqnarray*}
\sup_{1\leq i\leq n}
 \esp_{n}
  \biggl[
  \mathds{1}
  \left(X_i \in \mathcal{X}_0\right)
 (\widehat{m}_{in} - m(X_i))^4
 \biggr]
 &=&
 O_{\prob}
 \left(b_0^4+\frac{1}{nb_0^d}\right)^2,
 \\
 \sup_{1\leq i\leq n}
 \esp_{n}
 \biggl[
 \mathds{1}
 \left(X_i \in \mathcal{X}_0\right)
 (\widehat{m}_{in} - m(X_i))^6
 \biggr]
 &=&
 O_{\prob}
 \left(b_0^4+\frac{1}{nb_0^d}\right)^3.
\end{eqnarray*}
\label{BoundEspmchap}
\end{lem}

\begin{lem}
Assume that  $(A_4)$ and $(A_6)$ hold. Then, for any $1 \leq i
\neq j \leq n$, and for any $e$ in $\Rit$,
$$
\left(\widehat{m}_{in} - m(X_i), \varepsilon_i\right)
 \mbox{\it and }
 \left( \widehat{m}_{jn} - m(X_j), \varepsilon_j\right)
$$
are independent given $X_1, \ldots, X_n$,
 provided that $\| X_i - X_j \| \geq C b_0$,
for some constant $C>0$.
\label{Indep}
\end{lem}

\begin{lem}
Let $\Var_n(\cdot)$ and $Cov_n(\cdot)$ be respectively the
conditional variance and the conditional covariance given
$X_1,\ldots,X_n$, and set
\begin{eqnarray*}
\zeta_{in}
 =
\mathds{1} \left(X_i \in \mathcal{X}_0\right) (\widehat{m}_{in}
-m(X_i))^2 K_1^{(2)} \left( \frac{\varepsilon_i-e}{b_1} \right).
\end{eqnarray*}
Then under $(A_1)-(A_{9})$, we have, for $n$ going to infinity,
\begin{eqnarray*}
\sum_{i=1}^n
\Var_n\left( \zeta_{in}\right)
 &=&
O_{\prob}
\left(nb_1\right)
 \left(b_0^4+\frac{1}{nb_0^d}\right)^2,
 \\
 \sum_{i=1}^n
 \sum_{j=1\atop j\neq i}^n
  \Cov_n
 \left(
 \zeta_{in},\zeta_{jn}
 \right)
&=& O_{\prob}\left(n^2b_0^db_1^{7/2}\right)
\left(b_0^4+\frac{1}{nb_0^d}\right)^2.
\end{eqnarray*}
\label{sumzeta}
\end{lem}

\noindent All these lemmas are proved in Appendix A.

\subsection*{Proof of Theorem \ref{thm1}}
The proof of the theorem is based on the following equalities:
\begin{eqnarray}
\nonumber
 \widehat{f}_{n}(e)
 -
 \widetilde{f}_{n}(e)
&=& O_{\prob}
 \left[
 b_0^2
 +
 \left(
 \frac{1}{n}
 +
 \frac{1}{n^2b_0^db_1^3}
 \right)^{1/2}
 \right]
 +
 O_{\prob}
 \left[
 \frac{1}{(nb_1^5)^{1/2}}
 +
 \left(\frac{b_0^d}{b_1^3}\right)^{1/2}
 \right]
 \left(
 b_0^4
 +
 \frac{1}{nb_0^d}
 \right)
 \\
 &&
 +\;
 O_{\prob}
 \left[
 \frac{1}{b_1}
 +
 \left(\frac{b_0^d}{b_1^7}\right)^{1/2}
 \right]
 \left(b_0^4+\frac{1}{nb_0^d}\right)^{3/2},
 \label{fnchapfn1}
 \end{eqnarray}
 and
\begin{eqnarray}
\widetilde{f}_{n}(e)-f(e)
 =
 O_{\prob}
 \left(b_1^4+\frac{1}{nb_1}\right)^{1/2}.
 \label{fnchapfn2}
 \end{eqnarray}
Indeed, since $\widehat{f}_{n}(e)-f(e)
=(\widetilde{f}_{n}(e)-f(e))
+\widehat{f}_{n}(e)-\widetilde{f}_{n}(e)$, it then follows by
(\ref{fnchapfn2}) and (\ref{fnchapfn1}) that
\begin{eqnarray*}
\widehat{f}_{n}(e)-f(e)
 &=&
 O_{\prob}
 \left[
 b_1^4
 +
 \frac{1}{nb_1}
 +
 b_0^4
 +
 \frac{1}{n}
 +
 \frac{1}{n^2b_0^db_1^3}
 +
 \left(
 \frac{1}{(nb_1^5)^{1/2}}
 +
 \left(\frac{b_0^d}{b_1^3}\right)^{1/2}
 \right)^2
 \left(
 b_0^4
 +
 \frac{1}{nb_0^d}
 \right)^2
 \right]^{1/2}
 \\
 &&
 +\;
 O_{\prob}
 \left[
 \left(
 \frac{1}{b_1}
 +
 \left(\frac{b_0^d}{b_1^7}\right)^{1/2}
 \right)^2
 \left(b_0^4+\frac{1}{nb_0^d}\right)^3
 \right]^{1/2}.
 \end{eqnarray*}
This yields the result of the Theorem,  since  under ${\rm (A_8)}$
and ${\rm (A_9)}$, we have
\begin{eqnarray*}
 \frac{1}{n}
  =
 O\left(
 \frac{1}{nb_1}
 \right),
 \quad
\frac{1}{n^2b_0^db_1^3}
=
 O\left(\frac{b_0^d}{b_1^3}\right)
\left(b_0^4+\frac{1}{nb_0^d}\right)^2.
 \end{eqnarray*}
 Hence, it remains
to prove equalities (\ref{fnchapfn1}) and (\ref{fnchapfn2}). For
this, define $S_n$, $R_n$ and $T_n$ as in the statement of Lemma
\ref{STR}. Since $\widehat{\varepsilon}_i-\varepsilon_i =
-\left(\widehat{m}_{in}-m(X_{i})\right)$ and that $K_1$ is three
times continuously differentiable under ${\rm (A_7)}$,  the
third-order Taylor expansion with integral remainder gives
\begin{eqnarray*}
 \widehat{f}_{1n} (e)
 -
  \widetilde{f}_{n}(e)
 &=&
\frac{1}{b_1\sum_{i=1}^n \mathds{1}(X_i \in\mathcal{X}_0)}
\sum_{i=1}^n \mathds{1}\left(X_i\in \mathcal{X}_0\right)
\left[
 K_1\left(
 \frac{\widehat{\varepsilon}_i-e}{b_1}
 \right)
 -
 K_1\left(
 \frac{\varepsilon_i-e}{b_1}
 \right)
 \right]
 \\
& = &
 -\frac{1}{b_1\sum_{i=1}^n\mathds{1}(X_i \in\mathcal{X}_0)}
\left(
\frac{S_n}{b_1}
-
\frac{T_n}{2b_1^2}
+
\frac{R_n}{2b_1^3}
\right).
\end{eqnarray*}
Therefore, since
$$
\sum_{i=1}^n
\mathds{1}
\left(X_i \in \mathcal{X}_0\right)
 =
 n\left(
 \prob
 \left( X \in \mathcal{X}_0\right)
 +
 o_{\prob}(1)
 \right),
$$
 by the Law of large numbers, Lemma \ref{STR} then gives
\begin{eqnarray*}
\lefteqn{ \widehat{f}_{n}(e) - \widetilde{f}_{n}(e)
 =
O_{\prob}\left(\frac{1}{nb_1^2}\right)S_n
+
 O_{\prob}\left(\frac{1}{nb_1^3}\right) T_n
 +
 O_{\prob} \left(\frac{1}{nb_1^4}\right) R_n
 }
\\
&=&
 O_{\prob}
 \left[
 b_0^2
 \left(1+\frac{1}{(nb_1^3)^{1/2}}\right)
 +
 \left(
 \frac{1}{n}
 +
 \frac{1}{n^2b_0^db_1^3}
 \right)^{1/2}
 \right]
 \\
 &&
 +\;
 O_{\prob}
 \left[
 1
 +
 \frac{1}{(nb_1^5)^{1/2}}
 +
 \left(\frac{b_0^d}{b_1^3}\right)^{1/2}
 \right]
 \left(
 b_0^4
 +
 \frac{1}{nb_0^d}
 \right)
 +
 O_{\prob}
 \left[
 \frac{1}{b_1}
 +
 \left(\frac{b_0^d}{b_1^7}\right)^{1/2}
 \right]
 \left(b_0^4+\frac{1}{nb_0^d}\right)^{3/2}.
 \end{eqnarray*}
This yields (\ref{fnchapfn1}), since  under ${\rm (A_8)}$ and
${\rm (A_9)}$, we have $b_0\rightarrow 0$,
$nb_0^{d+2}\rightarrow\infty$ and $nb_1^3\rightarrow\infty$, so
that
\begin{eqnarray*}
b_0^2
 \left(1+\frac{1}{(nb_1^3)^{1/2}}\right)
&\asymp&
 O\left(b_0^2\right),
\quad
 \left(
 b_0^4
 +
 \frac{1}{nb_0^d}
 \right)
 =
 O\left(b_0^2\right),
\\
\left[
 1
 +
 \frac{1}{(nb_1^5)^{1/2}}
 +
 \left(\frac{b_0^d}{b_1^3}\right)^{1/2}
 \right]
 \left(
 b_0^4
 +
 \frac{1}{nb_0^d}
 \right)
 &=&
 O\left(b_0^2\right)
 +
 \left[
 \frac{1}{(nb_1^5)^{1/2}}
 +
 \left(\frac{b_0^d}{b_1^3}\right)^{1/2}
 \right]
 \left(
 b_0^4
 +
 \frac{1}{nb_0^d}
 \right).
\end{eqnarray*}
 For (\ref{fnchapfn2}), note that
\begin{eqnarray}
\esp_n \left[ \left( \widetilde{f}_{n}(e) -
 f(e)
\right)^2 \right] = \Var_n \left( \widetilde{f}_{n}(e) \right) +
\biggl( \esp_n \left[ \widetilde{f}_{n}(e) \right] - f(e)
\biggr)^2, \label{Quadbias}
\end{eqnarray}
with, using ${\rm (A_4)}$,
\begin{eqnarray*}
\Var_n \left( \widetilde{f}_{n}(e) \right)
 =
\frac{1}{
\left(
b_1\sum_{i=1}^n
\mathds{1}
\left(X_i\in\mathcal{X}_0\right)
\right)^2}
\sum_{i=1}^n
\mathds{1}
\left(X_i\in\mathcal{X}_0\right)
\Var
\left[
 K_1
 \left(
 \frac{\varepsilon-e}{b_1}
 \right)
 \right].
\end{eqnarray*}
Therefore, since the Cauchy-Schwarz inequality gives
\begin{eqnarray*}
\Var
\left[
 K_1\left(\frac{\varepsilon-e}{b_1}\right)
 \right]
\leq
 \esp
 \left[
 K_1^2\left(\frac{\varepsilon-e}{b_1}\right)
 \right]
 =
 b_1\int
 K_1^2(v)f(e+b_1v)
 dv,
\end{eqnarray*}
this bound and the equality above yield, under ${\rm (A_5)}$ and
${\rm (A_7)}$,
\begin{eqnarray}
 \Var_n\left(\widetilde{f}_{n}(e)\right)
 \leq
 \frac{C}{b_1\sum_{i=1}^n
 \mathds{1}\left(X_i\in\mathcal{X}_0\right)}
 =
 O_{\prob}\left(\frac{1}{nb_1}\right).
 \label{Quadbias1}
\end{eqnarray}
For the second term in (\ref{Quadbias}),  we have
\begin{eqnarray}
\esp_n\left[\widetilde{f}_{n}(e)\right] = \frac{1}{b_1\sum_{i=1}^n
 \mathds{1}\left(X_i\in\mathcal{X}_0\right)}
\sum_{i=1}^n \mathds{1}
\left(X_i\in\mathcal{X}_0\right)
\esp
\left[
 K_1\left(\frac{\varepsilon-e}{b_1}\right)
 \right].
 \label{Quadbias2}
\end{eqnarray}
By ${\rm (A_7)}$,  $K_1$  is symmetric, has  a compact support,
with
 $\int\! v K_1(v)=0$  and $\int \! K_1(v)dv=1$.
 Therefore, since $f(\cdot)$ has bounded continuous second
order derivatives under $(A_5)$, this yields for some
$\theta=\theta(e,b_1v)$,
\begin{eqnarray}
\nonumber
 \lefteqn{
 \esp
 \left[
 K_1\left(\frac{\varepsilon-e}{b_1}\right)
 \right]
 =
  b_1
\int
 K_1(v)
 f(e+b_1v)
 dv
 }
\\\nonumber
&=&
 b_1
 \int
 K_1(v)
\left[
 f(e)
+
 b_1 v f^{(1)}(e)
 +
 \frac{b_1^2v^2}{2}
 f^{(2)}(e+\theta b_1v)
 \right]
 dv
 \\
  &=&
 b_1f(e)
 +
 \frac{b_1^3}{2}
 \int v^2 K_1(v)
 f^{(2)}(e+\theta b_1v)
 dv.
 \label{MomK1}
\end{eqnarray}
Hence  this equality and (\ref{Quadbias2}) give
$$
\esp_n \left[\widetilde{f}_{n}(e)\right]
 =
 f(e)
 +
 \frac{b_1^2}{2}
 \int v^2 K_1(v)
 f^{(2)}(e+\theta b_1v)
 dv,
$$
so that
$$
\biggl( \esp_n \left[\widetilde{f}_{n}(e)\right] - f(e) \biggr)^2
 =
 O_{\prob}\left(b_1^4\right).
$$
Combining this result with  (\ref{Quadbias1}) and
(\ref{Quadbias}), we obtain, by the Tchebychev inequality,
$$
\widetilde{f}_{n}(e)-f(e) = O_{\prob}
\left(b_1^4+\frac{1}{nb_1}\right)^{1/2}.
$$
This proves (\ref{fnchapfn2}) and  achieves the proof of the
theorem.\eop

\subsection*{Proof of Theorem \ref{Optimbandw1}}
Recall that
\begin{eqnarray*}
 R_n(b_0, b_1)
=
 b_0^4
 +
 \left[
 \frac{1}{(nb_1^5)^{1/2}}
 +
 \left(\frac{b_0^d}{b_1^3}\right)^{1/2}
 \right]^2
 \left(
 b_0^4
 +
 \frac{1}{nb_0^d}
 \right)^2
 +
 \left[
 \frac{1}{b_1}
 +
 \left(\frac{b_0^d}{b_1^7}\right)^{1/2}
 \right]^2
 \left(b_0^4+\frac{1}{nb_0^d}\right)^3,
\end{eqnarray*}
and note that
$$
\left(\frac{1}{n^2b_1^3}\right)^{\frac{1}{d+4}}
=
\max
\left\lbrace
 \left(\frac{1}{n^2b_1^3}\right)^{\frac{1}{d+4}}
,
\left(\frac{1}{n^3b_1^7}\right)^{\frac{1}{2d+4}}
\right\rbrace
$$
if and only if $n^{4-d}b_1^{d+16}\rightarrow\infty$. To find the
order of $b_0^*$, we shall deal with the cases
$nb_0^{d+4}\rightarrow\infty$ and
 $nb_0^{d+4}=O(1)$.
\vskip 0.1cm\noindent First assume that
$nb_0^{d+4}\rightarrow\infty$. More precisely, we  suppose that
$b_0$ is in $\left[(u_n/n)^{1/(d+4)},+\infty\right)$, where $u_n
\rightarrow\infty$. Since $1/(nb_0^d) = O(b_0^4)$ for all these
$b_0$, we have
$$
\left(b_0^4+\frac{1}{nb_0^d}\right)^2
\asymp
\left(b_0^4\right)^2,
\quad \left(b_0^4+\frac{1}{nb_0^d}\right)^3
\asymp
 \left(b_0^4\right)^3.
$$
Hence the order of $b_0^*$ is computed  by minimizing the function
\begin{eqnarray*}
 b_0\rightarrow
  b_0^4
 +
 \left[
 \frac{1}{(nb_1^5)^{1/2}}
 +
 \left(\frac{b_0^d}{b_1^3}\right)^{1/2}
 \right]^2
 \left( b_0^4\right)^2
 +
 \left[
 \frac{1}{b_1}
 +
 \left(\frac{b_0^d}{b_1^7}\right)^{1/2}
 \right]^2
 \left(b_0^4\right)^3.
\end{eqnarray*}
Since this function is increasing with $b_0$, the minimum of $R_n
(\cdot,b_1)$ is achieved for $b_0^{**}=(u_n/n)^{1/(d+4)}$. We
shall prove later on that this choice of $b_0^{**}$ is irrelevant
compared to the one arising when $nb_0^{d+4}=O(1)$.

\vskip 0.1cm Consider now the case $nb_0^{d+4}=O(1)$ i.e
$b_0^4=O\left(1/(nb_0^d)\right)$. This  gives
\begin{eqnarray*}
 \left[
 \frac{1}{(nb_1^5)^{1/2}}
 +
 \left(\frac{b_0^d}{b_1^3}\right)^{1/2}
 \right]^2
 \left(
 b_0^4
 +
 \frac{1}{nb_0^d}
 \right)^2
 &\asymp&
  \left(
 \frac{1}{nb_1^5}
 +
 \frac{b_0^d}{b_1^3}
 \right)
 \left(
 \frac{1}{n^2b_0^{2d}}
 \right),
 \\
 \left[
 \frac{1}{b_1}
 +
 \left(\frac{b_0^d}{b_1^7}\right)^{1/2}
 \right]^2
 \left(b_0^4+\frac{1}{nb_0^d}\right)^3
 &\asymp&
 \left(
 \frac{1}{b_1^2}
 +
 \frac{b_0^d}{b_1^7}
 \right)
 \left(\frac{1}{n^3b_0^{3d}}\right).
 \end{eqnarray*}
Moreover if $nb_0^db_1^4\rightarrow\infty$, we have, since
$nb_0^{2d}\rightarrow\infty$ under ${\rm (A_8)}$,
$$
 \left(
 \frac{1}{nb_1^5}
 +
 \frac{b_0^d}{b_1^3}
 \right)
 \left(
 \frac{1}{n^2b_0^{2d}}
 \right)
 \asymp
 \frac{b_0^d}{b_1^3}
 \left(
 \frac{1}{n^2b_0^{2d}}
 \right),
 \quad
 \left(
 \frac{1}{b_1^2}
 +
 \frac{b_0^d}{b_1^7}
 \right)
 \left(\frac{1}{n^3b_0^{3d}}\right)
 =
 O\left(
 \frac{b_0^d}{b_1^3}
 \right)
 \left(
 \frac{1}{n^2b_0^{2d}}
 \right).
 $$
Hence the order of $b_0^*$ is obtained by finding the minimum of
the function $b_0^4+\left(1/n^2b_0^db_1^3\right)$. The
 minimization of this function  gives  a solution $b_0$ such that
$$
b_0
\asymp
\left(\frac{1}{n^2b_1^3}\right)^{\frac{1}{d+4}},
\quad
R_n(b_0,b_1)
\asymp
\left(\frac{1}{n^2b_1^3}\right)^{\frac{4}{d+4}}.
$$
This value  satisfies the constraints $nb_0^{d+4}=O(1)$ and
$nb_0^db_1^4\rightarrow\infty$ when
$n^{4-d}b_1^{d+16}\rightarrow\infty$.

\vskip 0.1cm\noindent If now $nb_0^{d+4}=O(1)$ but
$nb_0^db_1^4=O(1)$, we have, since $nb_0^{2d}\rightarrow\infty$,
\begin{eqnarray*}
\frac{1}{nb_1^5}
 \left(
\frac{1}{n^2b_0^{2d}} \right)
  =
 O\left(
 \frac{b_0^d}{b_1^7}
 \right)
 \left(\frac{1}{n^3b_0^{3d}}\right),
 \quad
 \frac{1}{b_1^2}
 \left(\frac{1}{n^3b_0^{3d}}\right)
 =
O\left(
 \frac{b_0^d}{b_1^3}
 \right)
 \left(
\frac{1}{n^2b_0^{2d}} \right)
 =
 O\left(
 \frac{b_0^d}{b_1^7}
 \right)
 \left(\frac{1}{n^3b_0^{3d}}\right).
\end{eqnarray*}
 In this case, $b_0^*$ is obtained by minimizing the function
$b_0^4+\left(1/n^3b_0^{2d}b_1^7\right)$, for which the solution
 $b_0$ verifies
$$
b_0
\asymp
\left(\frac{1}{n^3b_1^7}\right)^{\frac{1}{2d+4}},
\quad
R_n(b_0,b_1)
\asymp
\left(\frac{1}{n^3b_1^7}\right)^{\frac{4}{2d+4}}.
$$
This solution fulfills  the constraint $nb_0^db_1^4=O(1)$ when
$n^{4-d}b_1^{d+16}=O(1)$. Hence we can conclude that for
$b_0^4=O\left(1/(nb_0^d)\right)$, the bandwidth $b_0^*$ satisfies
$$
b_0^*
\asymp
\max
\left\lbrace
\left(\frac{1}{n^2b_1^3}\right)^{\frac{1}{d+4}}
,
\left(\frac{1}{n^3b_1^7}\right)^{\frac{1}{2d+4}}
\right\rbrace,
$$
which leads to
$$
R_n\left(b_0^*, b_1\right)
\asymp
\max
\left\lbrace
\left(\frac{1}{n^2b_1^3}\right)^{\frac{4}{d+4}}
,
\left(\frac{1}{n^3b_1^7}\right)^{\frac{4}{2d+4}}
 \right\rbrace.
$$
We need now to compare the solution $b_0^*$ to the candidate
 $b_0^{**}=(u_n/n)^{1/(d+4)}$ obtained when $nb_0^{d+4}\rightarrow\infty$.
  For this, we must  do a comparison between the
 orders of $R_n(b_0^*,b_1)$  and $R_n(b_0^{**},b_1)$. Since
 $R_n(b_0,b_1)\geq b_0^4$ for all $b_0$, we have
 $R_n(b_0^{**},b_1)\geq(u_n/n)^{4/(d+4)}$, so that
\begin{eqnarray*}
\frac{R_n(b_0^*, b_1)} {R_n(b_0^{**},b_1)}
&\leq&
 \left[
\left(\frac{1}{n^2b_1^3}\right)^{\frac{1}{d+4}}
 +
\left(\frac{1}{n^3b_1^7}\right)^{\frac{4}{2d+4}}
\right]
\left(\frac{n}{u_n}\right)^{\frac{4}{d+4}}
\\
&=&
o(1)
+
O\left(\frac{1}{u_n}\right)^{\frac{4}{d+4}}
\left(
\frac{1}{nb_1^{\frac{7(d+4)}{d+8}}}
\right)^{\frac{4(d+8)}{(2d+4)(d+4)}}
 = o(1),
\end{eqnarray*}
using the fact $n^{(d+8)}b_1^{7(d+4)}\rightarrow\infty$ by ${\rm
(A_9)}$ and  that $u_n\rightarrow\infty$.  This shows that
$R_n(b_0^*,b_1)\leq R_n(b_0^{**},b_1)$ for all $b_1$ and $n$ large
enough.
 Hence the theorem is proved, since $b_0^*$ is the best candidate for
 the minimization of $R_n(\cdot, b_1)$. \eop

\subsection*{Proof of Theorem \ref{Optimbandw2}}
Recall that  Theorem \ref{Optimbandw1} gives
\begin{eqnarray*}
AMSE(b_1)+R_n(b_0^*, b_1)
\asymp
r_1(b_1)+r_2(b_1)+r_3(b_1)
=
F(b_1),
\end{eqnarray*}
where
\begin{eqnarray*}
r_1(h)
&=&
h^4 + \frac{1}{nh},
\quad
\arg\min r_1(h)
\asymp
n^{-1/5}=h_1^*,
\quad
\min r_1(h)
\asymp
(h_1^*)^4=n^{-4/5},
\\
r_2(h)
&=&
h^4
+
\frac{1}{n^{\frac{8}{d+4}} h^{\frac{12}{d+4}}}
,
\quad
\arg\min r_2(h)
\asymp
n^{-\frac{2}{d+7}}
=
h_2^*,
\quad
\min r_2(h)
\asymp(h_2^*)^4
=
n^{-\frac{8}{d+7}},
\\
r_3(h)
&=&
h^4
+
\frac{1}{n^{\frac{12}{2d+4}}
h^{\frac{28}{2d+4}}},
\quad
\arg\min
r_3(h)
\asymp
n^{-\frac{3}{2d+11}}
=
h_3^*,
\quad
\min r_3(h)
\asymp(h_3^*)^4
=
n^{-\frac{12}{2d+11}}.
\end{eqnarray*}
Each $r_j(h)$ decreases on $\left[0,\arg\min r_j(h)\right]$ and
increases on $\left(\arg\min r_j(h),\infty\right)$ and that
$r_j(h)\asymp h^4$ on $\left(\arg\min r_j(h),\infty\right)$.
Moreover $\min r_2(h)=o\left(r_3(h)\right)$ and
$h_2^*=o\left(h_3^*\right)$ for all possible dimension $d$, so
that $\min\{r_2(h)+r_3(h)\}\asymp(h_3^*)^4=n^{-\frac{12}{2d+11}}$
and $\arg\min\{r_2(h)+r_3(h)\}\asymp h_3^*=n^{-\frac{3}{2d+11}}$.

\vskip 0.1cm Observe now that $\min\{r_2(h)+r_3(h)\}=O\left(\min
r_1(h)\right)$ is equivalent to
$n^{-\frac{12}{2d+11}}=O\left(n^{-4/5}\right)$ which holds if and
only if $d\leq 2$. Hence assume that $d\leq 2$. Since
$n^{-\frac{12}{2d+11}}=O\left(n^{-4/5}\right)$ also gives
$\arg\min\{r_2(h)+r_3(h)\}\asymp h_3^*=O\left(h_1^*\right)$, we
have
$$
\min F(b_1) \asymp n^{-4/5}
\;\;{\rm and}\;
\arg\min F(b_1)
\asymp
n^{-1/5}.
$$
The case $d>2$ is symmetric with
$$
\min F(b_1)
\asymp
n^{-\frac{12}{2d+11}}
\;\;{\rm and}\;
\arg\min F(b_1)
\asymp
n^{-\frac{3}{2d+11}}.
$$
This ends the proof of the theorem. \eop

\subsection*{Proof of Theorem \ref{normalite}}
Observe that the Tchebychev inequality gives
$$
\sum_{i=1}^n
\mathds{1}
 \left(X_i\in\mathcal{X}_0\right)
 =
n\prob\left(X\in\mathcal{X}_0\right)
 \left[
  1 + O_{\prob}
\left(\frac{1}{\sqrt{n}}\right)
\right],
$$
so that
$$
\widetilde{f}_{n}(e) = \left[ 1 + O_{\prob}
\left(\frac{1}{\sqrt{n}}\right) \right] f_n(e),
$$
where
$$
f_n(e)
 =
\frac{1}{nb_1\prob\left(X\in\mathcal{X}_0\right)}
\sum_{i=1}^n
\mathds{1}
\left(X_i\in\mathcal{X}_0\right)
 K_1\left(\frac{\varepsilon_i-e}{b_1}\right).
$$
Therefore
\begin{eqnarray}
\widehat{f}_{n}(e) - \esp f_n(e) = \left( f_n(e) - \esp f_n(e)
\right) + \left( \widehat{f}_{n}(e) - \widetilde{f}_{n}(e) \right)
+ O_{\prob} \left(\frac{1}{\sqrt{n}}\right)
 f_n(e).
\label{develop}
\end{eqnarray}
Let now $f_{in}(e)$ be as in the statement of Lemma \ref{Momfin1},
and note that $f_n(e) =(1/n)\sum_{i=1}^nf_{in}(e)$. The second and
the third claims of Lemma \ref{Momfin1} yield, since $nb_1$
diverges under ${\rm (A_9)}$,
\begin{eqnarray*}
\frac {
\sum_{i=1}^n
\esp\left| f_{in}(e)
-
\esp
f_{in}(e)
\right|^3}
{
\left(\sum_{i=1}^n
 \Var f_{in}(e)
\right)^{3/2}} \leq \frac{ \frac{C nf(e)}
{\prob\left(X\in\mathcal{X}_0\right)^2b_1^2}
 \displaystyle{\int}
\left| K_1(v) \right|^3 dv
+
 o\left( \frac{n}{b_1^2} \right)}
 {
\left( \frac{nf(e)}
{ \prob\left( X\in\mathcal{X}_0\right)b_1}
\displaystyle{\int}K_1^2(v) dv
+
 o\left(\frac{n}{b_1}\right)
\right)^{3/2}
 }
 =O(nb_1)^{-1/2}= o(1).
\end{eqnarray*}
 Hence the Lyapounov Central Limit Theorem for triangular arrays
 (see e.g Billingsley 1968, Theorem 7.3)
  gives, since $nb_1$ diverges under ${\rm (A_9)}$,
$$
\frac{f_n(e) - \esp f_{n}(e)} {\sqrt{\Var f_{n}(e)}} =
 \frac{f_n(e)
 -\esp f_{n}(e)}
 {\sqrt{\frac{\Var f_{in}(e)}{n}}}
 \stackrel{d}{\rightarrow}
 N\left(0,1\right).
$$
This yields, using the second result of Lemma \ref{Momfin1},
\begin{eqnarray}
\sqrt{nb_1} \left( f_n(e) - \esp f_{n}(e) \right)
\stackrel{d}{\rightarrow} N\left( 0, \frac{f(e)}
{\prob\left(X\in\mathcal{X}_0\right)} \int K_1^2(v) dv \right).
\label{develop1}
\end{eqnarray}
Moreover, note that for $nb_0^db_1^3\rightarrow\infty$ and
$nb_0^{2d}\rightarrow\infty$, we have
$$
 \frac{1}{nb_1^5}
 \left(
 \frac{1}{nb_0^d}
 \right)^2
 +
 \left(
 \frac{1}{b_1^2}
 +
 \frac{b_0^d}{b_1^7}
 \right)
 \left(\frac{1}{nb_0^d}\right)^3
=
O\left(
\frac{1}{n^2b_0^{d}b_1^3}
+
\frac{1}{n^3b_0^{2d}b_1^7}
\right).
$$
Therefore, since $b_0^4=O\left(1/(nb_0^d)\right)$,
$nb_0^db_1^3\rightarrow\infty$ and $nb_0^{2d}\rightarrow\infty$ by
$(\rm{A}_{10})$ and ${\rm (A_8)}$, the equality above and
(\ref{fnchapfn1}) ensure that
\begin{eqnarray*}
\widehat{f}_{n}(e) - \widetilde{f}_{n}(e)
 &=&
 O_{\prob}
 \left[
 b_0^4
 +
 \frac{1}{n}
 +
 \frac{1}{n^2b_0^db_1^3}
 +
 \left(
 \frac{1}{nb_1^5}
 +
 \frac{b_0^d}{b_1^3}
 \right)
 \left(
 \frac{1}{nb_0^d}
 \right)^2
 +
 \left(
 \frac{1}{b_1^2}
 +
 \frac{b_0^d}{b_1^7}
 \right)
 \left(\frac{1}{nb_0^d}\right)^3
 \right]^{1/2}
 \\
 &=&
 O_{\prob}
 \left(
 b_0^4
 +
 \frac{1}{n}
 +
 \frac{1}{n^2b_0^db_1^3}
 +
\frac{1}{n^3b_0^{2d}b_1^7}
 \right)^{1/2}.
\end{eqnarray*}
 Hence for $b_1$ going to $0$, we have
$$
\sqrt{nb_1} \left( \widehat{f}_{n}(e) - \widetilde{f}_{n}(e)
 \right)
 =
O_{\prob}
 \left[
 nb_1
 \left(
 b_0^4
 +
 \frac{1}{n}
 +
 \frac{1}{n^2b_0^db_1^3}
 +
\frac{1}{n^3b_0^{2d}b_1^7}
 \right)
\right]^{1/2}
 =
 o_{\prob}(1),
$$
since $nb_0^4b_1=o(1)$ and that $nb_0^db_1^3\rightarrow\infty$
under  $(\rm{A}_{10})$. Combining the above result with
(\ref{develop1}) and (\ref{develop}), we obtain
$$
\sqrt{nb_1}
 \left(
 \widehat{f}_{n}(e)
 -
 \esp f_n(e)
\right)
\stackrel{d}{\rightarrow}
 N\left(
 0,
\frac{f(e)} {\prob\left(X\in\mathcal{X}_0\right)} \int K_1^2(v) dv
 \right).
$$
This ends the proof the Theorem, since the first result of Lemma
\ref{Momfin1} gives
$$
\esp f_n(e) =
 \esp f_{1n}(e)
 =
  f(e) +
\frac{b_1^2}{2}f^{(2)}(e)
 \int v^2K_1(v) dv
  +
o\left(b_1^2\right) := \overline{f}_{n}(e). \eop
$$

 \setcounter{equation}{0} \setcounter{subsection}{0}
\setcounter{lem}{0}
\renewcommand{\theequation}{A.\arabic{equation}}
\renewcommand{\thesubsection}{A.\arabic{subsection}}
\begin{center}
\section*{Appendix A:  Proof of the  intermediate results}
\end{center}

\subsection*{Proof of Lemma \ref{Estig}}
First note that  by ${\rm (A_6)}$, we have
 $
 \int\!
z K_0(z) dz
 =0
 $
 and
 $
 \int\!
K_0(z) dz
 =1$.  Therefore, since $K_0$ is continuous and has a compact support,
${\rm (A_1)}$, ${\rm (A_2)}$ and the second-order Taylor expansion
yield, for  $b_0$ small enough and any $x$ in $\mathcal{X}_0$,
\begin{eqnarray*}
 \lefteqn{
 \left|
 \overline{g}_{n}(x)-g(x)
 \right|
  =
 \left|
 \frac{1}{b_0^d}
 \int
 K_0
 \left(\frac{z-x}{b_0}\right)
 g(z)
 dz
 -
 g(x)
 \right|
  =
 \left|
\int K_0(z) \left[g(x+b_0z)-g(x)\right]
 dz
 \right|
 }
 &&
\\
&=&
 \left|
 \int K_0(z)
 \left[
 b_0 g^{(1)}(x) z
 +
 \frac{b_0^2}{2}
 z g^{(2)} (x + \theta b_0 z)z^{\top}
\right] dz
 \right|,
 \;\theta = \theta (x,b_0 z)\in [0,1]
\\
&=& \left|
 b_0 g^{(1)}(x)
 \int
 z K_0(z) dz
 +
 \frac{b_0^2}{2}
  \int
z g^{(2)}(x + \theta b_0 z) z^{\top}
K_0(z)
dz
\right|
\\
& = & \frac{b_0^2}{2}
 \left|
 \int z g^{(2)}(x+\theta b_0z)
 z^{\top} K_0(z) dz
 \right|
 \leq C b_0^2.
\end{eqnarray*}
Therefore
$$
 \sup_{x\in\mathcal{X}_0}
 \left|
 \overline{g}_n (x) - g(x)
 \right|
 =
 O\left(b_0^2\right),
$$
which gives the first result of the lemma. For the two last
results of the lemma,  it is sufficient to show that
$$
 \sup_{x\in\mathcal{X}_0}
 \left|
 \widehat{g}_{n}(x)
  -
 \overline{g}_n(x)
 \right|
 =
 O_{\prob}
 \left(
  \frac{\ln n}{nb_0^d}
 \right)^{1/2},
$$
 since $\overline{g}_n (x)$ is asymptotically bounded away from
$0$ over
 $\mathcal{X}_0$
 and that
$|\overline{g}_n (x) - g(x) | = O (b_0^2) $ uniformly for $x$ in
$\mathcal{X}_0$. This follows from Theorem 1 in  Einmahl and Mason
(2005). \eop

\subsection*{Proof of Lemma \ref{Momfin1}}
For the first equality of the lemma, note that by ${\rm (A_4)}$,
${\rm (A_5)}$ and (\ref{MomK1}), we have
\begin{eqnarray*}
\esp \left[ f_{in}(e)\right] = \frac{1}{b_1} \esp \left[ K_1
\left( \frac{\varepsilon-e}{b_1} \right) \right] =
 f(e)
 +
 \frac{b_1^2}{2}
 \int v^2 K_1(v)
 f^{(2)}(e+\theta b_1v)
 dv.
\end{eqnarray*}
Therefore the Lebesgue Dominated Convergence Theorem gives, for
$b_1$ small enough,
\begin{eqnarray*}
\lefteqn{ \esp \left[ f_{in}(e) \right] - f(e) - \frac{b_1^2}{2}
f^{(2)}(e)
 \int
 v^2 K_1(v)
 dv
 }
 \\
 &=&
\frac{b_1^2}{2}
 \int
 v^2 K_1(v)
\left[ f^{(2)}\left(e+\theta b_1v\right) - f^{(2)}(e) \right]dv
\\
&=& o\left(b_1^2\right).
\end{eqnarray*}
This proves the first equality of the lemma. For the second and
third results of the lemma, the proofs are straightforward. Hence
they are omitted for the sake of brevity. \eop
%\begin{eqnarray*}
%\lefteqn{
%\Var [f_{in}(e)]
%=
%\esp \left[
%f_{in}^2(e)
%\right]
%-
%\left(\esp \left[ f_{in}(e)
%\right]
%\right)^2
%}
%\\
%&=&
%\frac{1}{b_1\prob\left(X\in\mathcal{X}_0\right)}
% \int
%f\left(e+b_1v\right)
%K_1^2(v)
%dv
%+
% O(1)
%\\
%&=&
%\frac{f(e)}
%{b_1\prob\left(X\in\mathcal{X}_0\right)}
%\int K_1^2(v) dv
%+
% o\left(\frac{1}{b_1}\right).
%\end{eqnarray*}
%
%The last statement of the lemma is immediate, since the Triangular
%and Convex inequalities give
%\begin{eqnarray*}
%\esp\left|
% f_{in}(e)
% -
% \esp f_{in}(e)
% \right|^3
%&\leq&
%C\esp\left| f_{in}(e)\right|^3
%\\
%&\leq&
%\frac {C f(e)}
%{b_1^2\prob^2\left(X\in\mathcal{X}_0\right) }
%\int
%\left| K_1(v)\right|^3
% dv
% +
% o\left(\frac{1}{b_1^2}\right).
% \eop
%\end{eqnarray*}

\subsection*{Proof of Lemma \ref{STR}}

The order of $S_n$ follows from Lemmas \ref{Betasum} and
\ref{Sigsum}. Indeed,  since
\begin{eqnarray*}
\mathds{1}(X_i \in \mathcal{X}_0)
 \left(
 \widehat{m}_{in} - m(X_i)
\right)
 &=&
 \frac{\mathds{1}(X_i \in \mathcal{X}_0)}
 {nb_0^d\widehat{g}_{in}}
 \sum_{j=1,j\neq i}^n
 \left( m(X_j)+ \varepsilon_j - m(X_i)\right)
 K_0\left(\frac{X_j-X_i}{b_0}\right)
 \\
 & =&
 \beta_{in} + \Sigma_{in},
\end{eqnarray*}
Lemmas \ref{Betasum} and \ref{Sigsum} imply that
\begin{eqnarray*}
S_n
 =
 O_{\prob}
\left[
b_0^2
\left(nb_1^2+(nb_1)^{1/2}\right)
 +
\left( nb_1^4 + \frac{b_1}{b_0^d} \right)^{1/2}
 \right],
\end{eqnarray*}
which gives the desired result for $S_n$.

 For the term $T_n$, the order is obtained by computing the
conditional mean and the conditional variance given
$X_1,\ldots,X_n$. To this end,  define for any $1\leq i\leq n$,
$$
\esp_{in} [\cdot]
=
\esp_n
\left[
 X_1,\ldots,X_n,\varepsilon_k,
 k\neq i
 \right].
$$
Therefore, since $(\widehat{m}_{in} - m(X_i))$ depends only upon
$\left(X_1,\ldots,X_n,\varepsilon_k,k\neq i\right)$, we have
\begin{eqnarray*}
 \esp_n [T_{n}]
 &=&
 \esp_{n}
 \left[
 \sum_{i=1}^n
 \esp_{in}
 \left[
 \mathds{1}
\left( X_i \in \mathcal{X}_0 \right)
 (\widehat{m}_{in} - m(X_i))^2
  K_1^{(2)}
 \left(\frac{\varepsilon_{i}-e}{b_1}\right)
 \right]
 \right]
 \\
 &=&
 \esp_{n}
  \left[
  \sum_{i=1}^n
  \mathds{1}
 \left( X_i \in \mathcal{X}_0\right)
(\widehat{m}_{in} - m(X_i))^2
 \esp_{in}
 \left[
 K_1^{(2)}
 \left(\frac{\varepsilon_{i}-e}{b_1}\right)
 \right]
 \right],
\end{eqnarray*}
with, using ${\rm (A_4)}$ and Lemma
\ref{MomderK}-(\ref{MomderK2}),
\begin{eqnarray*}
 \left|
\esp_{in}
 \left[
 K_1^{(2)}
 \left(\frac{\varepsilon_{i}-e}{b_1}\right)
 \right]
 \right|
 =
 \left|
 \int
  K_1^{(2)}
 \left(\frac{e-e}{b_1}
 \right)
 f(e)
 de
 \right|
\leq
 Cb_1^3.
\end{eqnarray*}
 Hence the equality above, the Cauchy-Schwarz inequality
  and  Lemma \ref{BoundEspmchap} yield
\begin{eqnarray}
\nonumber
 \left|
 \esp_{n}\left[T_{n}\right]
 \right|
&\leq&
 Cb_1^3
 \sum_{i=1}^n
  \esp_{n}
  \biggl[
  \mathds{1}\left( X_i \in \mathcal{X}_0\right)
 (\widehat{m}_{in} - m(X_i))^2
 \biggr]
 \\\nonumber
 &\leq&
 C n b_1^3
 \left(
 \sup_{1\leq i\leq n}
 \esp_{n}
  \biggl[
  \mathds{1}\left( X_i \in \mathcal{X}_0\right)
 (\widehat{m}_{in} - m(X_i))^4
 \biggr]
 \right)^{1/2}
 \\
 &\leq&
 O_{\prob}
 \left(nb_1^3\right)
 \left(b_0^4+\frac{1}{nb_0^d}\right).
 \label{Boundmean}
\end{eqnarray}
For the conditional variance of $T_n$,  Lemma \ref{sumzeta} gives
\begin{eqnarray*}
\Var_n (T_{n})
 &=&
 \sum_{i=1}^n
 \Var_n\left( \zeta_{in}\right)
 +
 \sum_{i=1}^n
 \sum_{j=1\atop j\neq i}^n
\Cov_n \left( \zeta_{in} , \zeta_{jn} \right)
\\
&=&
O_{\prob}\left(nb_1\right)
\left(b_0^4+\frac{b_1}{nb_0^d}\right)^2
+
O_{\prob}
\left(n^2b_0^db_1^{7/2}\right)
\left(b_0^4+\frac{1}{nb_0^d}\right)^2.
\end{eqnarray*}
Therefore, since $b_1$ goes to $0$,  the order above and
(\ref{Boundmean}) yield,
 applying the Tchebychev inequality,
\begin{eqnarray*}
T_{n}
 &=&
O_{\prob}
\left[
\left(nb_1^3\right)
\left(b_0^4+\frac{1}{nb_0^d}\right)
 +
\left(nb_1\right)^{1/2}
\left(b_0^4+\frac{b_1}{nb_0^d}\right)
+
\left(n^2b_0^db_1^{7/2}\right)^{1/2}
\left(b_0^4+\frac{1}{nb_0^d}\right)
\right]
\\
&=&
O_{\prob}
\left[
\left(
nb_1^3
+
\left(nb_1\right)^{1/2}
 +
\left(n^2b_0^db_1^{3}\right)^{1/2}
\right)
\left(b_0^4+\frac{1}{nb_0^d}\right)
\right].
\end{eqnarray*}
which gives the result for $T_n$.

 We now compute the order of $R_n$. For this, define
\begin{eqnarray*}
I_{in} &= & \int_{0}^{1}
 (1-t)^2
 K_1^{(3)}
 \left(
 \frac{
 \varepsilon_i-t(\widehat{m}_{in} -m(X_i))-e}{b_1}
 \right)
 dt,
 \\
R_{in}
& =&
 \mathds{1}
 \left(X_i\in\mathcal{X}_0\right)
 \left(\widehat{m}_{in}-m(X_i)\right)^3
 I_{in},
\end{eqnarray*}
and note that $R_n=\sum_{i=1}^n R_{in}$. The order of $R_n$ is
computed in a similar way as for $T_n$. Write
\begin{eqnarray*}
 \esp_n [R_{n}]
&=&
 \esp_n
 \left[
 \sum_{i=1}^n
 \esp_{in}\left[R_{in}\right]
\right]
\\
 &=&
 \esp_{n}
 \left[
 \sum_{i=1}^n
 \mathds{1}
 \left( X_i \in \mathcal{X}_0\right)
 (\widehat{m}_{in} - m(X_i))^3
 \esp_{in}
 \left[ I_{in}\right]\right],
\end{eqnarray*}
with, using ${\rm (A_4)}$ and Lemma
\ref{MomderK}-(\ref{MomderK3}),
\begin{eqnarray*}
 \left|
 \esp_{in}
 \left[ I_{in}\right]
 \right|
& =&
 \left|
 \int_{0}^{1}
 (1-t)^2
 \left[
 \int
  K_1^{(3)}
 \left(
 \frac{e-t(\widehat{m}_{in} -m(X_i))-e}{b_1}
 \right)f(e)
 de
 \right]
 dt
 \right|
 \\
&\leq&
 Cb_1^3.
\end{eqnarray*}
Therefore the Holder inequality and Lemma \ref{BoundEspmchap}
yield
\begin{eqnarray}
\nonumber
 \left|
 \esp_n\left[R_n\right]
 \right|
 &\leq&
 Cb_1^3
 \sum_{i=1}^n
 \esp_n
 \left[
\left|
\mathds{1}
\left(X_i\in\mathcal{X}_0\right)
\left(\widehat{m}_{in}-m(X_i)\right)
\right|^3
\right]
\\\nonumber
&\leq&
Cb_1^3
\sum_{i=1}^n
\esp_n^{3/4}
\left[
 \mathds{1}
\left(X_i\in\mathcal{X}_0\right)
\left(\widehat{m}_{in}-m(X_i)\right)^4
\right]
\\
&\leq&
 O_{\prob}
 \left(nb_1^3\right)
 \left(b_0^4+\frac{1}{nb_0^d}\right)^{3/2}.
 \label{Rn1}
\end{eqnarray}
For the conditional covariance of $R_n$, Lemma \ref{Indep} ensures
that
\begin{eqnarray}
\Var_n\left(R_n\right)
=
 \sum_{i=1}^n
 \Var_n\left(R_{in}\right)
 +
\sum_{i=1}^n
\sum_{j=1\atop j\neq i}^n
\biggl(\|X_i-X_j\|\leq
Cb_0\biggr)
\Cov_n\left(R_{in}, R_{jn}\right).
\label{Rn2}
\end{eqnarray}
Considering the  first term above, write
\begin{eqnarray*}
\Var_n\left(R_{in}\right)
 \leq
\esp_n\left[R_{in}^2\right] \leq
 \esp_{n}
 \biggl[
 \mathds{1}
 \left( X_i \in \mathcal{X}_0\right)
 (\widehat{m}_{in} - m(X_i))^6
 \esp_{in}
 \left[ I_{in}^2\right]
 \biggr],
\end{eqnarray*}
with, using  ${\rm (A_4)}$, the Cauchy-Schwarz inequality and
Lemma \ref{MomderK}-(\ref{MomderK3}),
\begin{eqnarray*}
 \esp_{in}
 \left[ I_{in}^2\right]
 &\leq&
 C\esp_{in}
 \left[
\int_{0}^{1}
 K_1^{(3)}
 \left(
 \frac{\varepsilon_i-t(\widehat{m}_{in} -m(X_i))-e}{b_1}
 \right)^2
 dt
 \right]
\\
&\leq&
 C\int_{0}^{1}
 \left[
 \int
 K_1^{(3)}
 \left(
 \frac{e-t(\widehat{m}_{in} -m(X_i))-e}{b_1}
 \right)^2
 f(e)
 de
 \right]
 dt
 \\
&\leq&
 Cb_1.
\end{eqnarray*}
Therefore
\begin{eqnarray*}
 \Var_n\left(R_{in}\right)
\leq
 Cb_1
 \esp_{n}
 \left[
 \mathds{1}
 \left( X_i \in \mathcal{X}_0\right)
 (\widehat{m}_{in} - m(X_i))^6
 \right],
 \end{eqnarray*}
uniformly in $i$. Hence  Lemma \ref{BoundEspmchap} imply that
\begin{eqnarray}
\nonumber
 \sum_{i=1}^n
 \Var_n\left(R_{in}\right)
 &\leq&
 Cnb_1
 \sup_{1\leq i\leq n}
 \esp_{n}
 \left[
 \mathds{1}
 \left( X_i \in \mathcal{X}_0\right)
 (\widehat{m}_{in} - m(X_i))^6
 \right]
 \\
&\leq&
 O_{\prob}
 \left(nb_1\right)
 \left(b_0^4+\frac{1}{nb_0^d}\right)^3.
 \label{Rn4}
\end{eqnarray}
For the second  term in (\ref{Rn2}),  we have
\begin{eqnarray*}
\left|
\Cov_n\left(R_{in}, R_{jn}\right)
 \right|
 &\leq&
 \left(
 \Var_n\left(R_{in}\right)
 \Var_n\left(R_{jn}\right)
 \right)^{1/2}
 \\
 &\leq&
 Cb_1
 \sup_{1\leq i\leq n}
 \esp_{n}
 \left[
 \mathds{1}
 \left( X_i \in \mathcal{X}_0\right)
 (\widehat{m}_{in} - m(X_i))^6
 \right].
\end{eqnarray*}
Hence from Lemma \ref{BoundEspmchap}  and the Tchebychev
inequality, we deduce
\begin{eqnarray*}
\lefteqn{
\sum_{i=1}^n
\sum_{j=1\atop j\neq i}^n
\biggl(\|X_i-X_j\|\leq Cb_0\biggr)
\left|
\Cov_n\left(R_{in}, R_{jn}\right)
\right|
}
\\
&\leq&
 O_{\prob}\left(b_1\right)
 \left(b_0^4+\frac{1}{nb_0^d}\right)^3
 \sum_{i=1}^n
 \sum_{j=1\atop j\neq i}^n
 \biggl(\|X_i-X_j\|\leq Cb_0\biggr)
 \\
 &\leq&
 O_{\prob}\left(b_1\right)
 \left(b_0^4+\frac{1}{nb_0^d}\right)^3
 \left(n^2b_0^d\right).
\end{eqnarray*}
This order, (\ref{Rn4}) and (\ref{Rn2}) give, since $nb_0^d$
diverges  under ${\rm (A_8)}$,
$$
\Var\left(R_n\right)
 =
 O_{\prob}
 \left(b_0^4+\frac{1}{nb_0^d}\right)^3
 \left(n^2b_0^db_1\right).
$$
Finally, with the help of this result and (\ref{Rn1}) we arrive at
\begin{eqnarray*}
R_n
&=&
O_{\prob}
\left[
\left(nb_1^3\right)
 \left(b_0^4+\frac{1}{nb_0^d}\right)^{3/2}
 +
\left(n^2b_0^db_1\right)^{1/2}
\left(b_0^4+\frac{1}{nb_0^d}\right)^{3/2}
\right]
\\
&=&
O_{\prob}
\left[
\left(
nb_1^3
+
\left(n^2b_0^db_1\right)^{1/2}
\right)
\left(b_0^4+\frac{1}{nb_0^d}\right)^{3/2}
\right].
 \eop
\end{eqnarray*}

\subsection*{Proof of Lemma \ref{MomderK}}
Set $h_p(e)=e^p f(e)$, $p\in[0,2]$. For  the first inequality of
(\ref{MomderK1}), note that  under ${\rm (A_5)}$ and ${\rm
(A_7)}$, the change of variable $e=e+b_1 v$ give, for any integer
$\ell\in[1, 3]$,
\begin{eqnarray}
\nonumber
 \left|
 \int
 K_1^{(\ell)}
 \left(\frac{e-e}{b_1}\right)^2
 e^pf(e) de
 \right|
 &=&
 \left|
 b_1
 \int
 K_1^{(\ell)}(v)^2 h_p(e+b_1v)
 dv
 \right|
 \\\nonumber
 &\leq&
 b_1
 \sup_{t\in\Rit}
 |h_p(t)|
 \int
 | K_1^{(\ell)}(v)^2|
 dv
 \\
 &\leq&
 Cb_1,
 \label{Ineg1}
\end{eqnarray}
which yields the first inequality in (\ref{MomderK1}). For the
second inequality in (\ref{MomderK1}), observe that $f(\cdot)$ has
a bounded  continuous derivative under ${\rm (A_5)}$, and that
$\int \! K_1^{(\ell)}(v)dv =0$ by ${\rm (A_7)}$. Therefore, since
$h_p(\cdot)$ has bounded second order derivatives under ${\rm
(A_6)}$, the Taylor inequality yields
\begin{eqnarray*}
\left| \int
 K_1^{(\ell)}
\left(\frac{e-e}{b_1}\right)
 e^pf(e)de
 \right|
 &=&
  b_1
 \left|
 \int
 K_1^{(\ell)}(v)
 \left[
 h_p(e+b_1v)-h_p(e)
 \right]
 \right|
  dv
 \\
 &\leq&
 b_1^2
 \sup_{t\in\Rit}|h_p^{(1)}(t)|
 \int
 |vK_1^{(\ell)}(v)|
  dv
 \leq
 Cb_1^2.
\end{eqnarray*}
which  proves (\ref{MomderK1}).
 The first inequalities of  (\ref{MomderK2}) and
 (\ref{MomderK3}) are given by (\ref{Ineg1}). The second bounds in
(\ref{MomderK2}) and (\ref{MomderK3}) are proved simultaneously.
For this, note that for any integer $\ell\in[2,3]$,
$$
\int
 K_1^{(\ell)}
 \left(\frac{e-e}{b_1}\right)
 h_p(e) de
 =
b_1 \int
 K_1^{(\ell)}(v)
h_p(e+b_1v) dv.
$$
Under ${\rm (A_7)}$,  $K_1(\cdot)$ is symmetric, has  a compact
support and two
 continuous derivatives, with
$\int \! K_1^{(\ell)}(v)dv=0$ and
 $\int\! v K_1^{(\ell)}(v) dv=0$.
 Hence the second order Taylor expansion applied to $h_p(\cdot)$ gives, for some
$\theta=\theta (e,b_1 v)\in[0,1]$,
\begin{eqnarray*}
 \lefteqn{
 \left|
 \int
 K_1^{(\ell)}
\left(
 \frac{e-e}{b_1}
 \right)
  h_p(e) de
 \right|
 =
 \left|
  b_1
 \int
 K_1^{(\ell)}(v)
\left[
 h_p(e+b_1v)
 -
  h_p(e)
 \right]
 dv
 \right|
 }
 &&
\\
&=&
\left|
b_1
\int
 K_1^{(\ell)}(v)
\left[
 b_1 v h_p^{(1)}(e)
 +
  \frac{b_1^2v^2}{2}
h_p^{(2)}(e+\theta b_1v) \right]
 dv
 \right|
 \\
  &=&
 \left|
\frac{b_1^3}{2}
\int
v^2 K_1^{(\ell)}(v)
 h_p^{(2)}(e+\theta b_1v)
 dv
\right|
\\
&\leq&
 \frac{b_1^3}{2}
 \sup_{t\in \Rit}|h_p^{(2)}(t)|
 \int
\left| v^2K_1^{(\ell)}(v)
 \right|
 dv
 \leq
 Cb_1^3,
\end{eqnarray*}
which completes the proof of the lemma.\eop

\subsection*{Proof of Lemma \ref{Betasum}}
By ${\rm (A_4)}$ and Lemma \ref{MomderK}-(\ref{MomderK1}) we have
\begin{eqnarray*}
\left|
\esp_n
\left[
 \sum_{i=1}^n
 \beta_{in}
 K_1^{(1)}
 \left(
\frac{\varepsilon_i-e}{b_1} \right)
 \right]
 \right|
 & = &
\left|
\esp
\left[
 K_1^{(1)}
  \left(
  \frac{\varepsilon -e}{b_1}
  \right)
  \right]
  \sum_{i=1}^n
  \beta_{in}
   \right|
  \leq C n b_1^2
 \max_{1 \leq i \leq n}
 \left| \beta_{in} \right|,
\\
\Var_n
 \left[
 \sum_{i=1}^n
  \beta_{in}
  K_1^{(1)}
  \left(
\frac{\varepsilon_i-e}{b_1}
 \right)
 \right]
&\leq&
  \sum_{i=1}^n
 \beta_{in}^2
 \esp
 \left[
 K_1^{(1)}
  \left(
   \frac{\varepsilon - e}{b_1}
\right)^2
 \right]
  \leq
 C n b_1
  \max_{1 \leq i \leq n}
  \left|
  \beta_{in}
\right|^2 .
\end{eqnarray*}
Hence the Tchebychev inequality gives
$$
\sum_{i=1}^n \beta_{in}
 K_1^{(1)}
 \left(
\frac{\varepsilon_i-e}{b_1} \right) =
 O_{\prob}
 \left( n
b_1^2 + (nb_1)^{1/2}
 \right)
 \max_{1 \leq i \leq n}
  \left|
\beta_{in}
 \right|,
$$
so that the lemma follows if we can prove that
\begin{equation}
\sup_{1 \leq i \leq n}
 \left|\beta_{in}\right|
 =
 O_{\prob}
\left( b_0^2\right),
 \label{BetasumTBP}
\end{equation}
as established now. For this, define
$$
\zeta_j (x)
 =
 \mathds{1}
 \left( x \in \mathcal{X}_0\right)
\left(m(X_j)-m(x)\right)
 K_0\left( \frac{X_j-x}{b_0}\right),
 \;\;
\nu_{in}(x)
 =
 \frac{1}{(n-1) b_0^d}
 \sum_{j=1, j\neq i}^n
 \left(
\zeta_j (x)-\esp[\zeta_j (x)]
 \right),
$$
and $\overline{\nu}_{n} (x)= \esp[\zeta_j(x)] / b_0^d $, so that
$$
\beta_{in} =
 \frac{n-1}{n} \frac{\nu_{in} (X_i)
 +
 \overline{\nu}_n (X_i)}{\widehat{g}_{in}}
\;.
$$
For  $\max_{1 \leq i \leq n} | \overline{\nu}_n (X_i) |$, first
observe that  a second-order Taylor expansion applied successively
to $g(\cdot)$ and $m(\cdot)$ give, for $b_0$ small enough, and for
any $x$, $z$ in $\mathcal{X}$,
\begin{eqnarray*}
 \lefteqn
 {
 \left[ m(x+b_0z)-m(x)\right]
 g(x+b_0z)
 }
\\
& =& \left[
 b_0 m^{(1)}(x) z
 +
 \frac{b_0^2}{2}
z m^{(2)}(x +\zeta_1 b_0 z)z^{\top}
 \right]
 \left[
 g(x)
 +
 b_0 g^{(1)}(x) z
 +
 \frac{b_0^2}{2}
 z g^{(2)}(x +\zeta_2 b_0z)z^{\top}
 \right],
  \end{eqnarray*}
for some $\zeta_1 =\zeta_1 (x,b_0 z)$ and $\zeta_2 =\zeta_2 (x,b_0
z)$ in $[0,1]$. Therefore, since $\int\!z
 K(z)dz=0$ under ${\rm (A_6)}$, it follows that,
 by ${\rm (A_1)}$, ${\rm (A_2)}$ and
 ${\rm (A_3)}$,
\begin{eqnarray}
\nonumber
 \max_{1 \leq i \leq n}
|\overline{\nu}_n (X_i)|
&\leq&
\sup_{x\in\mathcal{X}_0}
|\overline{\nu}_n (x)|
 =
 \sup_{x \in
 \mathcal{X}_0}
 \left|
 \int
 \left(m ( x + b_0 z) - m(x)\right)
 K_0 (z) g(x+b_0z)
 dz
\right|
\\
&\leq&
 Cb_0^2.
 \label{Betasum1}
\end{eqnarray}
 Consider now the term  $\max_{1\leq i \leq n}|\nu_{in}(X_i)|$.
  Using the Bernstein inequality (see e.g. Serfling (2002)), we have for any $t>0$,
\begin{eqnarray*}
\prob
\left(
\max_{1 \leq i \leq n}
| \nu_{in} (X_i)|
 \geq t
\right)
&\leq &
\sum_{i=1}^n \prob
 \left(
 | \nu_{in} (X_i) |
 \geq
t \right)
 \leq
 \sum_{i=1}^n
 \int
 \prob
 \left(
  | \nu_{in} (x) |
\geq t
 \left| X_i = x \right.
 \right) g (x)
  dx
\\
& \leq &
 2n \exp
 \left(
  - \frac{ (n-1) t^2 }
  { 2\sup_{x \in\mathcal{X}_0}
\Var (\zeta_j (x)/b_0^d) + \frac{4M}{3b_0^d} t}
 \right),
\end{eqnarray*}
where $M$ is such that $\sup_{x \in \mathcal{X}_0}|\zeta_j (x)|
\leq M$. Hence  ${\rm (A_2)}$, ${\rm (A_3)}$, ${\rm (A_6)}$ and
the standard Taylor expansion yield, for $b_0$ small enough,
$$
\sup_{x \in \mathcal{X}_0}
 | \zeta_j (x) |
 \leq C b_0,
 \;\;\;
\sup_{x \in \mathcal{X}_0}
 \Var (\zeta_j (x)/b_0^d)
\leq \frac{1}{b_0^d}
 \sup_{x \in \mathcal{X}_0}
 \int
 \left( m(x +b_0 z) - m(x) \right)^2
 K_0^2 (z) g(x+b_0z) dz
  \leq
  \frac{C
b_0^2}{b_0^d}\;,
$$
so that, for any $t \geq 0$,
$$
\prob
 \left(
 \max_{1 \leq i \leq n}
 | \nu_{in} (X_i) |
 \geq t
\right) \leq 2n
 \exp
 \left(
 - \frac{(n-1) b_0^d t^2 /b_0^2}{C + C
t/b_0}
\right).
$$
This gives
$$
\prob
\left(
\max_{1 \leq i \leq n}
|\nu_{in} (X_i)|
\geq
\left(
\frac{b_0^2 \ln n}{ (n-1) b_0^d}
\right)^{1/2}
t\right)
\leq
2n
\exp
\left(
 - \frac{ t^2 \ln n }
 {C + C t \left( \frac{\ln n}{
(n-1) b_0^d} \right)^{1/2} }
\right)
=
o(1),
$$
provided that $t$ is large enough and under ${\rm (A_8)}$. It then
follows that
$$
\max_{1 \leq i \leq n}
| \nu_{in} (X_i) |
 =
 O_{\prob}
 \left(
 \frac{b_0^2 \ln n}{ n b_0^d} \right)^{1/2}.
$$
This bound, (\ref{Betasum1}) and Lemma \ref{Estig}  show that
(\ref{BetasumTBP}) is proved,  since $b_0^2\ln
n/(nb_0^d)=O\left(b_0^4\right)$ under ${\rm (A_8)}$, and that
$$
 \beta_{in}
  =
 \frac{n-1}{n} \frac{\nu_{in} (X_i)
 +
 \overline{\nu}_n (X_i)}{\widehat{g}_{in}}\;.
 \eop
$$

\subsection*{Proof of Lemma \ref{Sigsum}}
Note that ${\rm (A_4)}$ gives that  $\Sigma_{in}$ is independent
of $\varepsilon_i$, and that $\esp_n[\Sigma_{in}]=0$. This yields
\begin{eqnarray}
 \esp_n
 \left[
 \sum_{i=1}^n
 \Sigma_{in}
 K_1^{(1)}
 \left(
 \frac{\varepsilon_i-e}{b_1}
 \right)
 \right]
  = 0.
\label{EspSigmai}
\end{eqnarray}
Moreover,  write
\begin{eqnarray}
\nonumber
 \lefteqn{
 \Var_n
  \left[
 \sum_{i=1}^n
 \Sigma_{in}
 K_1^{(1)}
  \left(
\frac{\varepsilon_i-e}{b_1}
 \right)
 \right]
 }
\\\nonumber
&=&
 \sum_{i=1}^n
 \Var_n
  \left[
 \Sigma_{in}
 K_1^{(1)}
 \left(
 \frac{\varepsilon_i-e}{b_1}
 \right)
 \right]
 +
 \sum_{i=1}^n
 \sum_{j=1\atop j\neq i}^n
 \Cov_n
 \left[
 \Sigma_{in}
 K_1^{(1)}
 \left(
 \frac{\varepsilon_{i}-e}{b_1}
 \right)
 ,
 \Sigma_{jn}
 K_1^{(1)}
 \left(
 \frac{\varepsilon_{j}-e}{b_1}
 \right)
 \right].
 \\
\label{VarSigm}
\end{eqnarray}
For the sum of variances in above, Lemma
\ref{MomderK}-(\ref{MomderK1}) and ${\rm (A_4)}$ give
\begin{eqnarray}
\nonumber \sum_{i=1}^n
 \Var_n
 \left[
 \Sigma_{in} K_1^{(1)}
 \left(
 \frac{\varepsilon_i-e}{b_1}
 \right)
 \right]
 &\leq&
 \sum_{i=1}^n
 \esp_n
\left[ \Sigma_{in}^2
 \right]
 \esp
 \left[
 K_1^{(1)}
  \left(
\frac{\varepsilon_i-e}{b_1} \right)^2
 \right]
\\\nonumber
&\leq&
 \frac{C b_1\sigma^2}{(nb_0^d)^2}
 \sum_{i=1}^n
 \sum_{j=1\atop j \neq i}^{n}
 \frac{\mathds{1}
 (X_i \in\mathcal{X}_0)}
{\widehat{g}_{in}^2}
 K_0^2
\left(\frac{X_j-X_i}{b_0}\right)
\\
 &\leq&
 \frac{C b_1\sigma^2}{nb_0^d}
 \sum_{i=1}^n
 \frac{\mathds{1}
 (X_i\in \mathcal{X}_0)
 \widetilde{g}_{in}}
{\widehat{g}_{in}^2}\;,
 \label{VarSigmai}
\end{eqnarray}
where $\sigma^2=\esp[\varepsilon^2]$ and
$$
\widetilde{g}_{in}
 =
 \frac{1}{n b_0^d}
 \sum_{j=1,j \neq i}^{n}
 K_0^2\left( \frac{X_j-X_i}{b_0}\right).
$$
For the sum of  conditional covariances in (\ref{VarSigm}), note
that
\begin{eqnarray*}
 \lefteqn{
 \sum_{i=1}^n
 \sum_{j=1\atop j\neq i}^n
\Cov_n \left[
 \Sigma_{in}
 K_1^{(1)}
 \left(
 \frac{\varepsilon_{i}-e}{b_1}
 \right)
 ,
 \Sigma_{jn}
 K_1^{(1)}
 \left(
  \frac{\varepsilon_{j}-e}{b_1}
 \right)
\right]
 }
\\
&=&
 \sum_{i=1}^n
 \sum_{j=1\atop j\neq i}^n
\esp_n
 \left[
 \Sigma_{in}
 \Sigma_{jn}
 K_1^{(1)}
 \left(
\frac{\varepsilon_{i}-e}{b_1}
 \right)
 K_1^{(1)}
  \left(
\frac{\varepsilon_{j}-e}{b_1}
 \right)
  \right]
\\
& = &
 \sum_{i=1}^n
 \sum_{j=1\atop j\neq i}^n
 \frac{\mathds{1}(X_{i}\in\mathcal{X}_0)
\mathds{1}(X_{j}\in\mathcal{X}_0)}
 {(n b_0^d)^2 \widehat{g}_{in} \widehat{g}_{jn}}
 \sum_{k=1\atop k\neq i}^n
\sum_{\ell=1\atop \ell\neq j}^n
 K_0
 \left(
\frac{X_{k}-X_{i}}{b_0}
 \right)
 K_0
 \left(
\frac{X_{\ell}-X_{j}}{b_0} \right) \esp
 \left[
 \xi_{ki}
 \xi_{\ell j}
 \right],
\end{eqnarray*}
where
$$
\xi_{ki}
 =
 \varepsilon_k
 K_1^{(1)}
  \left(
\frac{\varepsilon_{i}-e}{b_1}
 \right).
$$
Further, under ${\rm (A_4)}$, it is seen that for $k\neq\ell$,
$\esp[\xi_{ki}\xi_{\ell j}] =0 $ when $\Card\{i, j, k,
\ell\}\geq3$. Hence
 the symmetry of $K_0(\cdot)$ assumed in (A7) imply that
\begin{eqnarray*}
 \lefteqn{
 \sum_{i=1}^n
 \sum_{j=1\atop j\neq i}^n
\Cov_n \left[ \Sigma_{in}
 K_1^{(1)}
  \left(
\frac{\varepsilon_{i}-e}{b_1}
 \right),
 \Sigma_{jn}
K_1^{(1)}
 \left(
 \frac{\varepsilon_{j}-e}{b_1}
 \right)
 \right]
 }
 &&
\\
& = &
\sum_{i=1}^n
 \sum_{j=1\atop j\neq i}^n
\frac{\mathds{1}(X_{i} \in \mathcal{X}_0) \mathds{1} (X_{j}
 \in
\mathcal{X}_0)}
 {(n b_0^d)^2 \widehat{g}_{in} \widehat{g}_{jn}}
K_0^2 \left( \frac{X_{j}-X_{i}}{b_0} \right)
 \esp^2
 \left[
\varepsilon K_1^{(1)} \left( \frac{\varepsilon-e}{b_1} \right)
\right]
\\
&& + \sum_{i=1}^n
 \sum_{j=1\atop j\neq i}^n
 \frac{ \mathds{1} (X_{i}
\in \mathcal{X}_0) \mathds{1} (X_{j} \in \mathcal{X}_0)}
 {(nb_0^d)^2\widehat{g}_{in} \widehat{g}_{jn}}
 \sum_{k=1\atop k\neq  i, j}^n
 K_0 \left(\frac{X_{k}-X_{i}}{b_0}\right)
 K_0\left(\frac{X_{k}-X_{j}}{b_0}\right)
 \esp[\varepsilon^2]
\esp^2 \left[ K_1^{(1)}
 \left( \frac{\varepsilon-e}{b_1}
\right) \right].
\\
 %\label{CovSigmai}
\end{eqnarray*}
Therefore, since
$$
\sup_{1\leq i\leq n}
 \left(
\frac{\mathds{1}\left(X_i\in\mathcal{X}_0\right)}
{|\widehat{g}_{in}|} \right)
 =O_{\prob}(1),
$$
  by Lemma \ref{Estig},
  then Lemma \ref{MomderK}-(\ref{MomderK1})  gives
\begin{eqnarray}
\nonumber
 \lefteqn{
  \left|
 \sum_{i=1}^n
 \sum_{j=1\atop j\neq i}^n
\Cov_n
 \left[
  \Sigma_{in}
 K_1^{(1)}
 \left(
\frac{\varepsilon_{i}-e}{b_1}
 \right)
 ,
  \Sigma_{jn}
 K_1^{(1)}
 \left(
 \frac{\varepsilon_{j}-e}{b_1}
\right)
\right]
\right|
 }
\\
  &=&
 O_{\prob}
 \left(
 \frac{b_1^4}{n b_0^d}
 \right)
 \sum_{i=1}^n
 \mathds{1}(X_i\in
 \mathcal{X}_0)
 \widetilde{g}_{in}
 +
 O_{\prob} (b_1^4)
 \sum_{i=1}^n
 \mathds{1}(X_i\in \mathcal{X}_0)
 \widetilde{\widetilde{g}}_{in},
 \label{CovSigmaib}
\end{eqnarray}
where  $\widetilde{g}_{in}$ is defined as  in (\ref{VarSigmai})
and
$$
 \widetilde{\widetilde{g}}_{in}
 =
 \frac{1}{(n b_0^d)^2}
 \sum_{j=1\atop j\neq i}^n
 \sum_{k=1\atop k\neq j, i}^n
 K_0\left(\frac{X_k-X_i}{b_0}\right)
 K_0\left(\frac{X_k-X_j}{b_0}\right).
$$
In a completely similar way as done for Lemma \ref{Estig}, it can
be shown that $\widetilde{g}_{in}=O_{\prob}(1)$ uniformly in $i$
and for $n$ large enough. Therefore
 \begin{eqnarray}
 \sum_{i=1}^n
 \mathds{1}(X_i\in \mathcal{X}_0)
 \widetilde{g}_{in}
 =
 O_{\prob}(n).
 \label{CovSigmaib1}
 \end{eqnarray}
For the second term in (\ref{CovSigmaib}), the changes of
variables $x_1=x_3+b_0z_1$ and $x_2=x_3+b_0z_2$ give
\begin{eqnarray*}
 \esp
 \left[
 \sum_{i=1}^n
 \mathds{1}
\left(X_i\in\mathcal{X}_0\right)
\widetilde{\widetilde{g}}_{in}
 \right]
&\leq& \frac{Cn^3}{(nb_0^d)^2}
\esp
\left[
 K_0\left(\frac{X_3-X_1}{b_0}\right)
K_0\left(\frac{X_3-X_2}{b_0}\right)
 \right]
\\
&=&
\frac{Cn^3}{(nb_0^d)^2}
\int_{\mathcal{X}_0^3}
K_0\left(\frac{x_3-x_1}{b_0}\right)
K_0\left(\frac{x_3-x_2}{b_0}\right)
\prod_{k=1}^3
g(x_k)dx_k
 \\
 &\leq&
 \frac{Cn^3b_0^{2d}}{(nb_0^d)^2}=Cn,
\end{eqnarray*}
so that
$$
 \sum_{i=1}^n
 \mathds{1}(X_i\in \mathcal{X}_0)
 \widetilde{\widetilde{g}}_{in}
 =
 O_{\prob}(n).
$$
Hence from (\ref{VarSigm})-(\ref{CovSigmaib1}),
 we deduce
\begin{eqnarray*}
\lefteqn{
 \Var_n
 \left[
 \sum_{i=1}^n
 \Sigma_{in}
 K_1^{(1)}
 \left(
 \frac{\varepsilon_i-e}{b_1}
 \right)
 \right]
 }
 &&
 \\
 & =&
 O_{\prob}
 \left(
 \frac{b_1}{n b_0^d}
 \right)
 \sum_{i=1}^{n}
\mathds{1}(X_i\in \mathcal{X}_0)
\widetilde{g}_{in}
 +
 O_{\prob}
 \left(
 \frac{b_1^4}{n b_0^d}
 \right)
\sum_{i=1}^n
 \mathds{1}(X_i \in \mathcal{X}_0)
\widetilde{g}_{in}
 +
 O_{\prob}(b_1^4)
 \sum_{i=1}^n
\mathds{1}(X_i \in\mathcal{X}_0)
\widetilde{\widetilde{g}}_{in}
\\
&=&
 O_{\prob}
\left( \frac{b_1}{b_0^d}
 +
 \frac{b_1^4}{b_0^d}
 +
 nb_1^4
 \right)
 =
 O_{\prob}
 \left(
 \frac{b_1}{b_0^d}
 +
 nb_1^4
 \right).
\end{eqnarray*}
Finally, this order, (\ref{EspSigmai}) and the Tchebychev
inequality ensure that
$$
\sum_{i=1}^n
 \Sigma_{in}
 K_1^{(1)}
 \left(
 \frac{\varepsilon_i-e}{b_1}
 \right)
 =
 O_{\prob}
 \left(
 \frac{b_1}{b_0^d}+nb_1^4
 \right)^{1/2}. \eop
$$

\subsection*{Proof of Lemma \ref{BoundEspmchap}}

Define
\begin{eqnarray*}
g_{in}
 =
 \frac{1}{nb_0^d}
 \sum_{j=1, j\neq i}^n
 K_0^4\left(\frac{X_j-X_i}{b_0}\right),
 \quad
\widetilde{g}_{in}
 =
 \frac{1}{nb_0^d}
 \sum_{j=1, j\neq i}^n
 K_0^2\left(\frac{X_j-X_i}{b_0}\right).
\end{eqnarray*}
The proof of the lemma  is based on the following bound:
\begin{eqnarray}
\esp_n
 \biggl[
 \mathds{1}
 \left(X_i\in\mathcal{X}_0\right)
 \left(\widehat{m}_{in} - m(X_i)\right)^k
\biggr]
 \leq
 C
 \left[
 \beta_{in}^k
 +
 \frac{
 \mathds{1}
 \left(X_i\in\mathcal{X}_0\right)
 \widetilde{g}_{in}^{k/2}}
 {(nb_0^d)^{(k/2)}\widehat{g}_{in}^k}
 \right],
 \quad
 k\in\{4,6\}.
 \label{Espm}
 \end{eqnarray}
 Indeed, taking successively $k=4$ and $k=6$ in (\ref{Espm}),
 we have, by (\ref{BetasumTBP}), Lemma \ref{Estig} and
 ${\rm (A_8)}$,
\begin{eqnarray*}
 \sup_{1\leq i\leq n}
 \esp_n
 \biggl[
 \mathds{1}\left(X_i\in\mathcal{X}_0\right)
 \left(\widehat{m}_{in} - m(X_i)\right)^4
\biggr]
 &=&
 O_{\prob}
 \left(
 b_0^8
  +
 \frac{1}{(nb_0^d)^2}
 \right)
 =
 O_{\prob}
 \left(b_0^4+\frac{1}{nb_0^d}\right)^2,
 \\
 \sup_{1\leq i\leq n}
 \esp_n
 \biggl[
 \mathds{1}
 \left(X_i\in\mathcal{X}_0\right)
 \left(\widehat{m}_{in} - m(X_i)\right)^6
\biggr]
 &=&
 O_{\prob}
 \left(
 b_0^{12}
  +
 \frac{1}{(nb_0^d)^3}
 \right)
 =
 O_{\prob}
 \left(b_0^4+\frac{1}{nb_0^d}\right)^3,
\end{eqnarray*}
which gives the desired results of the lemma.
 Hence it remains to prove (\ref{Espm}).
 For this,  define $\beta_{in}$ and  $\Sigma_{in}$  respectively as in
 the statement of Lemmas \ref{Betasum} and \ref{Sigsum}.
 Since
 $\mathds{1}(X_i \in\mathcal{X}_0)
 \left(\widehat{m}_{in}- m(X_i)\right)
 =\beta_{in}+\Sigma_{in}$, and that
 $\beta_{in}$ depends only upon
 $\left(X_1,\ldots,X_n\right)$, this gives, for $k\in\{4,6\}$
\begin{eqnarray}
\esp_n
 \biggl[
 \mathds{1}
 (X_i \in \mathcal{X}_0)
 \left(\widehat{m}_{in}- m(X_i)\right)^k
 \biggr]
  \leq
 C\beta_{in}^k
 +
 C\esp_n\left[\Sigma_{in}^k\right].
 \label{Espm5}
\end{eqnarray}
The order  of the second term of bound (\ref{Espm5}) is computed
by applying Theorem 2 in Whittle (1960) or the
Marcinkiewicz-Zygmund inequality (see e.g Chow and Teicher, 2003,
p. 386). These inequalities show that for linear form
$L=\sum_{j=1}^n a_j\zeta_j$ with independent mean-zero random
variables
 $\zeta_1,\ldots,\zeta_n$, it holds that, for any $k\geq 1$,
 $$
 \esp
 \left|L^k\right|
 \leq
  C(k)
  \left[
 \sum_{j=1}^n
  a_j^2
\esp^{2/k} \left|\zeta_j^k \right|
 \right]^{k/2},
 $$
where $C(k)$ is a positive real depending only on $k$. Now,
observe that for any $i\in[1,n]$,
$$
\Sigma_{in}
 =
 \sum_{j=1, j\neq i}^n
 \sigma_{jin},
 \quad
 \sigma_{jin}
 =
 \frac{\mathds{1}
 \left( X_i \in \mathcal{X}_0\right)}
 {nb_0^d\widehat{g}_{in}}
 \varepsilon_j
 K_0\left(\frac{X_j-X_i}{b_0}\right).
$$
Since under ${\rm (A_4)}$, the $\sigma_{jin}$'s, $j\in[1,n]$, are
centered independent variables given $X_1,\ldots,X_n$,  this
yields, for any $k\in\{4,6\}$,
\begin{eqnarray*}
 \esp_n
 \left[
 \Sigma_{in}^k\right]
 \leq
 C\esp\left[\varepsilon^k\right]
 \left[
 \frac{\mathds{1}
 \left( X_i \in \mathcal{X}_0\right)}
 {(nb_0^d)^2\widehat{g}_{in}^2}
 \sum_{j=1}^n
  K_0^2
 \left(
 \frac{X_j-X_i}{b_0}
 \right)
 \right]^{k/2}
 \leq
 \frac{C\mathds{1}
 \left(X_i\in\mathcal{X}_0\right)\widetilde{g}_{in}^{k/2}}
 {(nb_0^d)^{(k/2)}\widehat{g}_{in}^k}\;.
\end{eqnarray*}
 Hence this  bound and
(\ref{Espm5})  give
$$
 \esp_n
 \biggl[
 \mathds{1}(X_i \in \mathcal{X}_0)
 \left(\widehat{m}_{in}- m(X_i)\right)^k
 \biggr]
 \leq C\left[
 \beta_{in}^k
  +
 \frac{
 \mathds{1}
 \left(X_i\in\mathcal{X}_0\right)
 \widetilde{g}_{in}^{k/2}}
 {(nb_0^d)^{(k/2)}\widehat{g}_{in}^k}
 \right],
$$
which proves (\ref{Espm}) and then completes the proof of the
lemma. \eop

\subsection*{Proof of Lemma \ref{Indep}}

Since $K_0(\cdot)$ has a compact support under ${\rm (A_6)}$,
there is a $C>0$ such that $\| X_i - X_j \| \geq C b_0$ implies
that for any integer number $k$ of $[1,n]$, $K_0 ( (X_k -
X_i)/b_0) = 0$ if $K_0 ( (X_j - X_k)/b_0) \neq 0$. Let $D_j
\subset [1,n]$ be such that an integer number $k$ of $[1,n]$ is in
$D_j$ if and only if $K_0 ( (X_j - X_k)/b_0) \neq 0$. Abbreviate
$\prob (\cdot| X_1, \ldots,X_n)$ into $\prob_n$ and assume that
$\| X_i - X_j \| \geq C b_0$ so that $D_i$ and $D_j$ have an empty
intersection. Note also that taking $C$ large enough ensures that
$i$ is not in $D_j$ and $j$ is not in $D_i$. It then follows,
under ${\rm (A_4)}$ and since $D_i$ and $D_j$ only depend upon
$X_1,\ldots,X_n$,
\begin{eqnarray*}
\lefteqn{
 \prob_n
 \biggl(
 \left(
 \widehat{m}_{in} - m(X_i),\varepsilon_i
 \right)
 \in A \mbox{ \rm and }
 \left(
\widehat{m}_{jn} - m (X_j),\varepsilon_j
\right)
 \in B \biggr)
 }
&&
\\
& = & \prob_n \left( \left( \frac{\sum_{k \in D_i \setminus\{i\}}
 \left( m(X_{k}) - m(X_i) + \varepsilon_{k} \right)
  K_0
\left( (X_{k} - X_i)/b_0 \right)} {\sum_{k \in D_i \setminus\{i\}}
 K_0 \left( (X_{k} - X_i)/b_0 \right)}
 ,
 \varepsilon_i\right) \in A \right.
\\
&& \;\;\;\;\;\;\;\;\;\;\;\;\;\;\;\;\;\;\;\; \left. \mbox{ \rm and}
\left( \frac{ \sum_{\ell \in D_j \setminus \{j \}}
 \left(
m(X_{\ell}) - m(X_j) + \varepsilon_{\ell} \right)
 K_0 \left(
(X_{\ell} - X_j)/b_0 \right) } { \sum_{\ell \in D_j \setminus \{j
\}} K_0 \left( (X_{\ell} - X_j)/b_0 \right)} ,
\varepsilon_j\right) \in B \right)
\\
& = & \prob_n \left( \left( \frac{ \sum_{k \in D_i \setminus \{i
\}} \left( m(X_{k}) - m(X_i) + \varepsilon_{k}\right)
 K_0 \left(
(X_{k} - X_i)/b_0 \right)} {\sum_{k \in D_i \setminus \{i \}}
K_0\left( (X_{k} - X_i)/b_0 \right)} ,
 \varepsilon_i \right) \in
A \right)
\\
&& \;\;\;\;\;\;\;\;\;\;\;\;\;\;\;\;\;\;\;\;
 \times\;
 \prob_n
 \left(
\left( \frac{ \sum_{\ell \in D_j \setminus \{j \}} \left(
m(X_{\ell}) - m(X_j) + \varepsilon_{\ell} \right) K_0 \left(
(X_{\ell} - X_j)/b_0 \right) } { \sum_{\ell \in D_j \setminus \{j
\}} K_0\left( (X_{\ell} - X_j)/b_0\right) } , \varepsilon_j
\right) \in B \right)
\\
& = & \prob_n
 \left(
 \left(\widehat{m}_{in} - m(X_i), \varepsilon_i \right)
 \in A
 \right)
\times \prob_n \left( \left(\widehat{m}_{jn} - m
(X_j),\varepsilon_j \right)
 \in B
\right).
\end{eqnarray*}
This gives the result of Lemma \ref{Indep}, since both
$\left(\widehat{m}_{in} - m (X_i), \varepsilon_i\right)$ and
$\left(\widehat{m}_{jn} - m (X_j), \varepsilon_j\right)$ are
independent given $X_1, \ldots, X_n$. \eop

\subsection*{Proof of Lemma \ref{sumzeta}}
Since $\widehat{m}_{in} - m(X_i)$ depends only upon
 $\left(X_1,\ldots,X_n,\varepsilon_k, k\neq i\right)$,
 we have
\begin{eqnarray*}
\sum_{i=1}^n
 \Var_n
 \left(\zeta_{in}\right)
  \leq
 \sum_{i=1}^n
 \esp_n
 \left[
 \zeta_{in}^2
 \right]
 =
 \sum_{i=1}^n
 \esp_n
 \left[
 \mathds{1}
 \left(X_i\in\mathcal{X}_0\right)
 \left(\widehat{m}_{in} - m(X_i)\right)^4
  \esp_{in}
  \left[
  K_1^{(2)}
  \left(\frac{\varepsilon_{i}-e}{b_1}\right)^2
  \right]
 \right],
 \end{eqnarray*}
with, using Lemma \ref{MomderK}-(\ref{MomderK2}),
\begin{eqnarray*}
 \esp_{in}
 \left[
 K_1^{(2)}
 \left(\frac{\varepsilon_{i}-e}{b_1}\right)^2
 \right]
 =
 \int
  K_1^{(2)}
 \left(
 \frac{e-e}{b_1}
 \right)^2
 f(e) de
 \leq
 C b_1.
\end{eqnarray*}
 Therefore these bounds and Lemma \ref{BoundEspmchap} give
\begin{eqnarray*}
\sum_{i=1}^n
 \Var_n\left(\zeta_{in}\right)
 &\leq&
 Cb_1
 \sum_{i=1}^n
 \esp_n
 \biggl[
 \mathds{1}
 \left(X_i\in\mathcal{X}_0\right)
 (\widehat{m}_{in} - m(X_i))^4
 \biggr]
 \\
 &\leq&
 Cnb_1
 \sup_{1\leq i\leq n}
 \esp_n
 \biggl[
 \mathds{1}
 \left(X_i\in\mathcal{X}_0\right)
 (\widehat{m}_{in} - m(X_i))^4
 \biggr]
 \\
 &\leq&
 O_{\prob}\left(nb_1\right)
 \left(b_0^4+\frac{1}{nb_0^d}\right)^2.
\end{eqnarray*}
which yields the desired result for the conditional variance.

 We now prepare to compute the order of the conditional covariance.
 Observe that Lemma \ref{Indep} gives
\begin{eqnarray*}
\sum_{i=1}^n \sum_{j=1\atop j\neq i}^n
\Cov_n \left(
\zeta_{in},\zeta_{jn}
 \right)
 =
  \sum_{i=1}^n
  \sum_{j=1\atop j\neq i}^n
  \mathds{1}
  \biggl(\left\|X_i - X_j\right\|<C b_0\biggr)
  \biggl(
  \esp_n
  \left[
  \zeta_{in}
  \zeta_{jn}
  \right]
  -
 \esp_n\left[\zeta_{in}\right]
 \esp_n\left[\zeta_{jn}\right]
 \biggr).
\end{eqnarray*}
 The order of the term above  is derived from the following
equalities:
\begin{eqnarray}
 \sum_{i=1}^n
 \sum_{j=1\atop j\neq i}^n
 \mathds{1}
 \biggl(
 \left\|X_i - X_j \right\|<C b_0
 \biggr)
 \esp_n\left[\zeta_{in}\right]
 \esp_n\left[\zeta_{jn}\right]
 &=&
 O_{\prob}
 \left(n^2b_0^db_1^6\right)
\left(b_0^4+ \frac{1}{nb_0^d}\right)^2,
  \label{Covzeta2}
  \\
  \sum_{i=1}^n
  \sum_{j=1\atop j\neq i}^n
  \mathds{1}
  \biggl(
  \left\|X_i - X_j\right\|<C b_0
  \biggr)
  \esp_n
  \left[
  \zeta_{in}
  \zeta_{jn}
  \right]
  &=&
  O_{\prob}
 \left(n^2b_0^db_1^{7/2}\right)
\left(b_0^4+ \frac{1}{nb_0^d}\right)^2.
 \label{Covzeta1}
\end{eqnarray}
 Indeed, since $b_1$ goes to $0$ under ${\rm (A_9)}$,
 (\ref{Covzeta2}) and (\ref{Covzeta1}) yield
\begin{eqnarray*}
 \sum_{i=1}^n
 \sum_{j=1\atop j\neq i}^n
 \Cov_n
 \left(\zeta_{in},\zeta_{jn}\right)
 &=&
O_{\prob}
 \left[
 \left(n^2b_0^db_1^6\right)
 \left(b_0^4+ \frac{1}{nb_0^d}\right)^2
 +
 \left(n^2b_0^db_1^{7/2}\right)
 \left(b_0^4+ \frac{1}{nb_0^d}\right)^2
\right]
\\
&=& O_{\prob}
 \left(n^2b_0^db_1^{7/2}\right)
 \left(b_0^4+ \frac{1}{nb_0^d}\right)^2,
\end{eqnarray*}
 which gives the result for the conditional
covariance. Hence,  it remains to prove (\ref{Covzeta2}) and
(\ref{Covzeta1}).  For (\ref{Covzeta2}), note that by ${\rm
(A_4)}$ and Lemma \ref{MomderK}-(\ref{MomderK2}), we have
 \begin{eqnarray*}
 \left|
 \esp_{n}
 \left[\zeta_{in}\right]
 \right|
 &=&
 \left|
 \esp_n
 \left[
 \mathds{1}
 \left(X_i\in\mathcal{X}_0\right)
 (\widehat{m}_{in} - m(X_i))^2
 \esp_{in}
 \left[
 K_1^{(2)}
 \left(
 \frac{\varepsilon_i-e}{b_1}
 \right)
 \right]
 \right]
 \right|
 \\
 &\leq&
 Cb_1^3
 \biggl(
  \esp_n
 \biggl[
 \mathds{1}
 \left(X_i\in\mathcal{X}_0\right)
 (\widehat{m}_{in} - m(X_i))^4
 \biggr]
 \biggr)^{1/2}.
 \end{eqnarray*}
Hence from this bound and Lemma \ref{BoundEspmchap} we deduce
\begin{eqnarray*}
\sup_{1\leq i, j\leq n}
 \left|
 \esp_{n}
 \left[\zeta_{in}\right]
 \esp_n
\left[\zeta_{jn}\right]
 \right|
 &\leq&
 Cb_1^6
 \sup_{1\leq i \leq n}
 \esp_n
 \biggl[
 \mathds{1}
 \left(X_i\in\mathcal{X}_0\right)
 (\widehat{m}_{in} - m(X_i))^4
 \biggr]
 \\
 &\leq&
 O_{\prob}
 \left(b_1^6\right)
\left(b_0^4+ \frac{1}{nb_0^d}\right)^2.
\end{eqnarray*}
Therefore, since
\begin{eqnarray}
\sum_{i=1}^n
 \sum_{j=1\atop j\neq i}^n
 \mathds{1}
 \biggl(
 \| X_i - X_j \|< C b_0
 \biggr)
 =
 O_{\prob}(n^2 b_0^d),
 \label{Markov}
 \end{eqnarray}
 by the Tchebychev inequality gives, it then follows that
\begin{eqnarray*}
  \sum_{i=1}^n
 \sum_{j=1\atop j\neq i}^n
 \mathds{1}
 \biggl(
 \| X_i - X_j \|< C b_0
 \biggr)
 \esp_{n}
 \left[\zeta_{in}\right]
 \esp_{n}
 \left[\zeta_{jn}\right]
  =
 O_{\prob}
 \left(n^2b_0^db_1^6\right)
\left(b_0^4+ \frac{1}{nb_0^d}\right)^2,
\end{eqnarray*}
which   proves (\ref{Covzeta2}).
 For (\ref{Covzeta1}), set $Z_{in}= \mathds{1}
\left(X_i\in\mathcal{X}_0\right)\left(\widehat{m}_{in} -
m(X_i)\right)^2$, and note that for $i\neq j$, we have
 \begin{eqnarray}
 \esp_n\left[\zeta_{in}\zeta_{jn}\right]
 =
 \esp_n
 \left[
 Z_{in}
  K_1^{(2)}
 \left(
 \frac{\varepsilon_j-e}{b_1}
 \right)
  \esp_{in}
 \left[
 Z_{jn}
  K_1^{(2)}
 \left(
 \frac{\varepsilon_i-e}{b_1}
 \right)
 \right]
 \right],
 \label{Prodzeta}
 \end{eqnarray}
where
\begin{eqnarray}
\nonumber
 \lefteqn{
 \esp_{in}
 \left[
 Z_{jn}
  K_1^{(2)}
 \left(
 \frac{\varepsilon_i-e}{b_1}
 \right)
 \right]
}
\\\nonumber
&=&
 \beta_{jn}^2
 \esp_{in}
 \left[
 K_1^{(2)}
 \left(
 \frac{\varepsilon_i-e}{b_1}
 \right)
 \right]
 +
 2\beta_{jn}
 \esp_{in}
 \left[
 \Sigma_{jn}
 K_1^{(2)}
 \left(
 \frac{\varepsilon_i-e}{b_1}
 \right)
 \right]
 +
 \esp_{in}
 \left[
 \Sigma_{jn}^2
 K_1^{(2)}
 \left(
 \frac{\varepsilon_i-e}{b_1}
 \right)
 \right].
 \\
\label{Covzeta3}
\end{eqnarray}
The first term of  (\ref{Covzeta3}) is treated by using Lemma
\ref{MomderK}-(\ref{MomderK2}). This gives
\begin{eqnarray}
\left| \beta_{jn}^2
 \esp_{in}
 \left[
 K_1^{(2)}
 \left(
 \frac{\varepsilon_i-e}{b_1}
 \right)
 \right]
 \right|
 \leq
 Cb_1^3
 \beta_{jn}^2.
 \label{Covzeta4}
\end{eqnarray}
Since  under ${\rm (A_4)}$, the $\varepsilon_j$'s are independent
centered variables, and  are independent of the  $X_j$'s, the
second term of (\ref{Covzeta3}) equals
\begin{eqnarray*}
\lefteqn{
 2\beta_{jn}\frac{\mathds{1}\left(X_j\in\mathcal{X}_0\right)}
 {nb_0^d\widehat{g}_{jn}}
\sum_{k=1, k\neq j}^n
 K_0
 \left(
 \frac{X_k-X_j}{b_0}
 \right)
\esp_{in}
 \left[
 \varepsilon_k
   K_1^{(2)}
 \left(
 \frac{\varepsilon_i-e}{b_1}
 \right)
 \right]
 }
 \\
 &=&
  2\beta_{jn}
 \frac{\mathds{1}\left(X_j\in\mathcal{X}_0\right)}
 {nb_0^d\widehat{g}_{jn}}
 K_0
 \left(
 \frac{X_i-X_j}{b_0}
 \right)
\esp_{in}
 \left[
 \varepsilon_i
   K_1^{(2)}
 \left(
 \frac{\varepsilon_i-e}{b_1}
 \right)
 \right].
\end{eqnarray*}
Therefore, since $K_0$ is bounded under ${\rm (A_6)}$, the
equality above and Lemma \ref{MomderK}-(\ref{MomderK2}) imply that
\begin{eqnarray}
\left|
2\beta_{jn}
\esp_{in}
 \left[
 \Sigma_{jn}
   K_1^{(2)}
 \left(
 \frac{\varepsilon_i-e}{b_1}
 \right)
 \right]
\right|
 \leq
 Cb_1^3
\left| \beta_{jn}\right|
\frac{\mathds{1}\left(X_j\in\mathcal{X}_0\right)}
 {nb_0^d\widehat{g}_{jn}}.
 \label{Covzeta5}
\end{eqnarray}
For the  last term of (\ref{Covzeta3}), we have
\begin{eqnarray*}
\lefteqn{
 \esp_{in}
 \left[
 \Sigma_{jn}^2(x)
   K_1^{(2)}
 \left(\frac{\varepsilon_i-e}{b_1}\right)
 \right]
}
\\
&=&
\frac{1}{(nb_0^d\widehat{g}_{jn})^2}
 \sum_{k=1\atop k\neq j}^n
\sum_{\ell=1\atop\ell\neq j}^n
 K_0\left(\frac{X_k-X_j}{b_0}\right)
 K_0\left(\frac{X_{\ell}-X_j}{b_0}\right)
 \esp_{in}
 \left[
 \varepsilon_k
 \varepsilon_{\ell}
  K_1^{(2)}
 \left(\frac{\varepsilon_i-e}{b_1}\right)
 \right]
 \\
 &=&
 \frac{1}{(nb_0^d\widehat{g}_{jn})^2}
 \sum_{k=1, k\neq j}^n
 K_0^2\left(\frac{X_k-X_j}{b_0}\right)
 \esp_{in}
 \left[
 \varepsilon_k^2
  K_1^{(2)}
 \left(\frac{\varepsilon_i-e}{b_1}\right)
 \right],
\end{eqnarray*}
with, using Lemma \ref{MomderK}-(\ref{MomderK2}),
\begin{eqnarray*}
\lefteqn
 {
 \left|
 \esp_{in}
 \left[
 \varepsilon_k^2
  K_1^{(2)}
 \left(\frac{\varepsilon_i-e}{b_1}\right)
 \right]
 \right|
}
\\
&\leq& \max
 \left\lbrace
\sup_{e\in\Rit}
\left|
 \esp_{in}
 \left[
 \varepsilon^2
  K_1^{(2)}
 \left(\frac{\varepsilon-e}{b_1}\right)
 \right]
 \right|,
 \;
 \esp[\varepsilon^2]
\sup_{e\in\Rit}
 \left|
 \esp_{in}
 \left[
 K_1^{(2)}
 \left(\frac{\varepsilon-e}{b_1}\right)
 \right]
 \right|
\right\rbrace
\\
&\leq&
 C b_1^3.
\end{eqnarray*}
Therefore
$$
 \left|
 \esp_{in}
 \left[
 \Sigma_{jn}^2
   K_1^{(2)}
 \left(
 \frac{\varepsilon_i-e}{b_1}
 \right)
 \right]
 \right|
 \leq
\frac{Cb_1^3}{(nb_0^d\widehat{g}_{jn})^2}
\sum_{k=1, k\neq j}^n
K_0^2\left(\frac{X_k-X_j}{b_0}\right).
$$
Substituting this bound, (\ref{Covzeta5}) and (\ref{Covzeta4}) in
(\ref{Covzeta3}), we obtain
$$
\left|
 \esp_{in}
 \left[
 Z_{jn}
  K_1^{(2)}
 \left(
 \frac{\varepsilon_i-e}{b_1}
 \right)
 \right]
\right|
 \leq
 Cb_1^3 M_n,
$$
where
$$
M_n
=
\sup_{1\leq j\leq n}
\left[
\beta_{jn}^2
+
\left|
\beta_{jn}\right|
\frac{\mathds{1}
\left(X_j\in\mathcal{X}_0\right)}
{nb_0^d\widehat{g}_{jn}}
+
\frac{1}{(nb_0^d\widehat{g}_{jn})^2}
\sum_{k=1, k\neq j}^n
K_0^2\left(\frac{X_k-X_j}{b_0}\right)
 \right].
$$
Hence from (\ref{Prodzeta}), the Cauchy-Schwarz inequality, Lemma
\ref{BoundEspmchap} and Lemma \ref{MomderK}-(\ref{MomderK2}), we
deduce
\begin{eqnarray*}
\lefteqn{
\sum_{i=1}^n
\sum_{j=1\atop j\neq i}^n
 \mathds{1}
 \biggl(
 \| X_i - X_j \|< C b_0
 \biggr)
\left|
 \esp_n
 \left[ \zeta_{in}\zeta_{jn} \right]
 \right|
 }
\\
&\leq&
C M_nb_1^3
\sum_{i=1}^n
\sum_{j=1\atop j\neq i}^n
 \mathds{1}
 \biggl(
 \| X_i - X_j \|< C b_0
 \biggr)
\esp_n
 \left|
 Z_{in} K_1^{(2)}
\left(\frac{\varepsilon_j-e}{b_1}\right) \right|
\\
&\leq&
CM_nb_1^3
 \sum_{i=1}^n
 \sum_{j=1\atop j\neq i}^n
 \mathds{1}
 \biggl(
 \| X_i - X_j \|< C b_0
 \biggr)
\esp_n^{1/2} \left[ Z_{in}^2 \right]
 \esp_n^{1/2}
 \left[
 K_1^{(2)}
 \left(
 \frac{\varepsilon_j-e}{b_1}
 \right)^2
 \right]
 \\
 &\leq&
 M_n b_1^3
 O_{\prob}
 \left(b_0^4+\frac{1}{nb_0^d}\right)
 (b_1)^{1/2}
\sum_{i=1}^n \sum_{j=1\atop j\neq i}^n
\biggl( \mathds{1}
\left(\|X_i-X_j\|\leq Cb_0\right)
\biggr).
\end{eqnarray*}
Further, using (\ref{BetasumTBP}) and Lemma \ref{Estig}, it can be
shown that
$$
M_n =
 O_{\prob}
 \left( b_0^4 + \frac{b_0^2}{nb_0^d}
  +
\frac{1}{nb_0^d}
 \right)
= O_{\prob}
\left(b_0^4+ \frac{1}{nb_0^d}\right).
$$
Therefore, substituting this order in the inequality above, and
using (\ref{Markov}), we arrive at
\begin{eqnarray*}
\sum_{i=1}^n \sum_{j=1\atop j\neq i}^n
 \mathds{1}
 \biggl(
 \| X_i - X_j \|< C b_0
 \biggr)
\esp_n
\left[\zeta_{in}\zeta_{jn}\right]
=
 O_{\prob}
\left(n^2b_0^db_1^{7/2}\right)
\left(b_0^4+\frac{1}{nb_0^d}\right)^2,
\end{eqnarray*}
which proves (\ref{Covzeta1}) and  completes the proof of the
lemma. \eop

\end{document}